\documentclass[a4paper, 12pt, final]{amsart}

\makeatletter
\numberwithin{equation}{section}
\numberwithin{figure}{section}
\usepackage[a4paper]{geometry}
\usepackage{amsmath,amsfonts,amssymb,amsthm,epsfig}

\usepackage{color,xcolor}
\usepackage[final,linkcolor = blue,citecolor = blue,colorlinks=true]{hyperref}
\usepackage[leqno]{amsmath}

\usepackage{amscd}
\usepackage{mathrsfs}
\usepackage{}
\usepackage[all]{xy}
\usepackage[OT1]{fontenc}
\usepackage[latin1]{inputenc}
\usepackage{amssymb,latexsym}
\usepackage{type1cm,ae}
\usepackage[notref,notcite]{showkeys}
\usepackage{graphicx}
\usepackage{xspace}
\usepackage{setspace}
\usepackage{enumerate}
\usepackage[english]{babel}
\usepackage{bbm}
\usepackage{esint}
\usepackage{soul}

\soulregister{\cite}7
\soulregister{\citep}7
\soulregister{\citet}7
\soulregister{\ref}7
\soulregister{\pageref}7
\sethlcolor{yellow}

\theoremstyle{definition}
\newtheorem{theorem}{Theorem}[section]
\newtheorem{proposition}[theorem]{Proposition}

\newtheorem{lemma}[theorem]{Lemma}

\numberwithin{equation}{section}

\newcommand{\e}{\varepsilon}

\newcommand*{\loc}{\mathrm{loc}}

\newcommand*{\Wert}{\mathord{\mbox{|\kern-1.5pt|\kern-1.5pt|}}}

\newcommand{\N}{\mathbb{N}}	
\newcommand{\R}{\mathbb{R}}	

\DeclareMathOperator{\diam}{diam}

\newcommand{\Div}{\mathrm{Div}}
\newcommand{\diag}{\mathrm{diag}}

\def\XXint#1#2#3{{\setbox0=\hbox{$#1{#2#3}{\int}$}
		\vcenter{\hbox{$#2#3$}}\kern-.5\wd0}}

\makeatother

\usepackage{babel}
\begin{document}
	
	\title[Regularity of a geometrically nonlinear system]{$L^p$-regularity of a geometrically nonlinear system in supercritical dimensions} 
	
	\author[C.-Y. Guo, M.-L. Liu and C.-L. Xiang]{Chang-Yu Guo, Ming-Lun Liu and Chang-Lin Xiang}
	
	\address[Chang-Yu Guo]{Research Center for Mathematics and Interdisciplinary Sciences, Shandong University, 266237, Qingdao, China and  Department of Physics and Mathematics, University of Eastern Finland, 80101, Joensuu, Finland}
	\email{guocybnu@gmail.com}
	
	\address[Ming-Lun Liu]{Research Center for Mathematics and Interdisciplinary Sciences, Shandong University 266237,  Qingdao, P. R. China and and  Frontiers Science Center for Nonlinear Expectations, Ministry of Education, P. R. China}
	\email{minglunliu2021@163.com}
	
	\address[Chang-Lin Xiang]{Three Gorges Mathematical Research Center, China Three Gorges University,  443002, Yichang,  P. R. China} \email{changlin.xiang@ctgu.edu.cn}
	

	\thanks{$^*$Corresponding author: Ming-Lun Liu}
	\thanks{C.-Y. Guo and M.-L. Liu are supported by the Young Scientist Program of the Ministry of Science and Technology of China (No.~2021YFA1002200), the NSFC grant (No.~12101362), the Taishan Scholar Project and the NSF of Shandong Province (No.~ZR2022YQ01). C.-L. Xiang is supported by the NFSC grant (No.~12271296) and  the NSF of Hubei province (No. 2024AFA061).}
	
	\date{}
	
	\begin{abstract} 
	In a recent work, Gastel and Neff introduced an interesting system from a geometrically nonlinear flat cosserat micropolar model and established interior regularity in the critical dimension. Inspired by their work on this flat Cosserat model, in this article, we establish both interior regularity and sharp $L^p$ regularity for their system in supercritical dimensions.   
	\end{abstract}

	\maketitle
	
\section{Introduction}\label{Introduction}
Motivated by applications to a geometrically nonlinear flat Cosserat shell model in continuum mechanics, Gastel and Neff~\cite{Gastel-Neff 2022} introduced the following interesting system: 
\begin{align}\label{detailed eq:G-N system 1st eq}
	&\Div S(\nabla m,R) =0,  \\
	\label{detailed eq:G-N system 2nd eq}
	&\Delta R-\Omega_{R}\cdot\nabla R-\mathrm{skew} \left(\nabla m\circ S(\nabla m,R)\right)R =0, 
\end{align}
where the unknown functions $(m,R)\in W^{1,2}(B^n,\R^N\times SO(N))$ with $\nabla m,\nabla R\in M^{2,n-2}(B^n)$, 
$$
\Omega_R=-R \nabla R^{\mathrm{T}}\in M^{2,n-2}(B^n, \R^n\otimes so(N)),
$$
$S:\mathbb{R}^{N\times n}\times \mathbb{R}^{N\times N} \to \mathbb{R}^{N\times n} $ is a mapping given in \eqref{eq:def for S} and $\mathrm{skew} \left(\nabla m\circ S(\nabla m,R)\right)$ is defined by \eqref{eq:def for skew m s}. Here, $B^n$ is the unit ball in $\R^n$ and $N\ge n$.

This system couples a harmonic map type equation \eqref{detailed eq:G-N system 2nd eq} with a uniformly elliptic equation \eqref{detailed eq:G-N system 1st eq}. In their main result, Gastel and Neff obtained the following interior regularity result in the planar case, that is, $n=2$ and $N=3$. 

\textbf{Theorem A} [{Theorem 1.1,~\cite{Gastel-Neff 2022}}]\label{thmA:Gastel-Neff}
Every weak solution $(m,R)\in W^{1,2}(B^2,\R^3\times SO(3))$ of \eqref{detailed eq:G-N system 1st eq}-\eqref{detailed eq:G-N system 2nd eq} (with $n=2, N=3$) is smooth.

As one easily observes, the function $\Omega_R$ from \eqref{detailed eq:G-N system 2nd eq} is the same as that in 
\begin{equation}\label{eq:harmonic map eq}
	\Delta R-\Omega_R\cdot \nabla R=0,
\end{equation}
which models harmonic mappings from $B^n$ to $SO(N)\subset \R^{N\times N}$. In the planar case, that is when $n=2$, the regularity of harmonic mappings into manifolds was first obtained by Morrey in his seminal work~\cite{Morrey-1948} on Plateau's problem in Riemannian manifolds. In particular, he showed that minimizing harmonic mappings are locally H\"older continuous and thus are smooth when the Riemannian metric is smooth. This regularity result was later extended to weakly harmonic mappings by Hel\'ein in his celebrated work; see~\cite{Helein-2002} for a comprehensive introduction on it. Partially building on H\'elein's idea, in another significant work, Riv\`ere~\cite{Riviere-2007} successfully rewrote the harmonic mapping equation \eqref{eq:harmonic map eq} as a conservation law, from which regularity and compactness follow routinely. An important discovery of Rivi\`ere~\cite{Riviere-2007} is that the specific form of $\Omega_R$ (as that of harmonic mappings) is not really essential. The algebraic anti-symmetry of $\Omega_R$ is sufficient for finding the conservation law, based on the earlier seminal work of Uhlenbeck~\cite{Uhlenbeck-1982}. We recommend the interested readers to~\cite{Riviere-2012} for a comprehensive exploration of Rivi\`ere's conservation law approach and to ~\cite{Gastel-2019,Li-Wang-2022-IMRN} for related works on regularity of the Cosserat shell model. 

Unfortunately, because of the appearance of an extra term $\mathrm{skew}\left(\nabla m\circ S(\nabla m,R)\right)R$ in \eqref{detailed eq:G-N system 2nd eq}, the powerful conservation law approach of Hel\'ein and Rivi\`ere does not work. The key observation by Gastel and Neff~\cite{Gastel-Neff 2022} is that this extra term is indeed the product of a gradient term $\nabla m$, a divergence-free vector field $S(\nabla m,R)$ and a bounded term $R$, which makes it enjoy better properties than merely being in $L^1$. Adapting the method of Rivi\`ere and Struwe~\cite{RS-2008}, Gastel and Neff succeeded in deriving local H\"older regularity and thus also smoothness of weak solutions to \eqref{detailed eq:G-N system 1st eq}-\eqref{detailed eq:G-N system 2nd eq}. 

As for the harmonic mapping equation \eqref{eq:harmonic map eq}, partial regularity results have been known since the fundamental work of Schoen and Uhlenebck~\cite{Shoen-Uhlenbeck-82}, where the important $\epsilon$-regularity was established for minimizing harmonic mappings into manifolds. This work was later extended to the case of stationary harmonic mappings into spheres by Evans~\cite{Evans-1991} and into general manifolds by Bethuel~\cite{Bethuel-1993}. An alternative (but slightly more general) proof was obtained later by Rivi\`ere and Struwe~\cite{RS-2008}, partially based on observations from~\cite{Riviere-2007}. It remains, however, open whether one can extend the conservation law approach of Rivi\`ere~\cite{Riviere-2007} to study (partial) regularity of harmonic mappings in supercritical dimensions; see~\cite{Guo-Xiang 2 order} for some partial progress in this direction. It is then natural to ask whether one can derive partial regularity for weak solutions of \eqref{detailed eq:G-N system 1st eq}-\eqref{detailed eq:G-N system 2nd eq} in supercritical dimensions. Indeed, in~\cite[Remark 6.8]{Gastel-Neff 2022}, Gastel and Neff formulated it as an open question:
\smallskip 

\textbf{Question B:} Can we establish a partial regularity theory for weak solutions of \eqref{detailed eq:G-N system 1st eq}-\eqref{detailed eq:G-N system 2nd eq} in supercritical dimensions? Can we find good estimates in the natural Morrey spaces? 
\medskip 

Note that in the case of (stationary) harmonic mappings, due to the monotonicity formula, the gradient of a harmonic map lies in the Morrey space $M^{2,n-2}$; see for instance~\cite{RS-2008}. Thus for Question B, the natural Morrey spaces for a partial regularity theory would be $\nabla m,\nabla R\in M^{2,n-2}$. 

Our main motivation for this paper is to provide an affirmative answer to Question B of Gastel and Neff. Our first main result concerns interior regularity of weak solutions. 
\begin{theorem}[Interior regularity]\label{thm:interior regularity}
	There exists some $\epsilon=\epsilon(n,N)>0$ such that if $(m,R)\in W^{1,2}(B^n,\R^N\times SO(N))$ is a weak solution of \eqref{detailed eq:G-N system 1st eq}-\eqref{detailed eq:G-N system 2nd eq} with
	$$
	\|\nabla m\|_{M^{2,n-2}(B^n)}\le\e\quad \text{and}\quad \|\nabla R\|_{M^{2,n-2}(B^n)}\le\e,
	$$
	then it is smooth in $B_{\frac{1}{2}}$.
\end{theorem}

As was pointed out in~\cite{Gastel-Neff 2022}, the main difficulty in establishing Theorem \ref{thm:interior regularity} is the (local) H\"older continuity of weak solutions. In the proof of Theorem \ref{thm:interior regularity}, we closely follow the idea of Gaste and Neff~\cite{Gastel-Neff 2022} and thus rely on the ideas of Rivi\`ere and Struwe~\cite{RS-2008}, except that we refine some of the estimates using  Hardy-BMO inequalities (see Lemma \ref{Hardy-BMO duality by L^p-L^q-Morrey} below). Once the H\"older continuity is established, smoothness follows routinely. 

Motivated by applications in the associated heat flow and energy identity \cite{Eells-Sampson-1964,Jost-Liu-Zhu-2019-MathAnn,Li-Zhu-2012-PJM,Struwe-1985,Wang-Wei-Zhang-2017}, the $L^p$-regularity theory for harmonic mapping type equations has attracted great interest in the past decades; see for instance~\cite{Lamm-Sharp-2016-CPDE,Schoen-84,Sharp-Topping-2013-TAMS,Sharp-2014,Moser-2015-TAMS} and the references therein.   
Motivated by these works and potential applications, we aim at obtaining optimal $L^p$-regularity estimates for the corresponding inhomogeneous system:
\begin{align}\label{detailed eq:G-N system 3rd eq}
	&\Div S(\nabla m,R) =0,  \\
	\label{detailed eq:G-N system 4th eq}
	&\Delta R-\Omega_{R}\cdot\nabla R-\mathrm{skew} \left(\nabla m\circ S(\nabla m,R)\right)R =f, 
\end{align}
where the inhomogeneous term $f\in L^p(B^n,\R^{N\times N})$. 

Our second main result provides optimal $L^p$-regularity estimates for weak solutions of \eqref{detailed eq:G-N system 3rd eq}-\eqref{detailed eq:G-N system 4th eq}.

\begin{theorem}[$L^p$-regularity]\label{thm:Lp regularity}
	Suppose $f\in L^p(B^n,\R^{N\times N})$ for some $p\in(\frac{n}{2},\infty)$. There exists some $\epsilon=\epsilon(n,N)>0$ such that if $(m,R)\in W^{1,2}(B^n,\R^N\times SO(N))$ is a weak solution of \eqref{detailed eq:G-N system 3rd eq}-\eqref{detailed eq:G-N system 4th eq} with
	$$
	\|\nabla m\|_{M^{2,n-2}(B^n)}\le\e\quad \text{and}\quad \|\nabla R\|_{M^{2,n-2}(B^n)}\le\e,
	$$ 
	then $(m,R)\in W^{2,p}(B_{\frac{1}{2}})$. Furthermore, there exists some $C=C(n,N,p)>0$ such that  
	\begin{equation}\label{eq:f R-m W^2,p estimate}
		\|m\|_{W^{2,p}(B_{\frac{1}{2}})}+\|R\|_{W^{2,p}(B_{\frac{1}{2}})}\le C\left(\|f\|_{L^p(B_1)}+1\right)^2.
	\end{equation} 
\end{theorem}

As an immediate consequence of Theorem \ref{thm:Lp regularity}, we know that if $f=0$, then the solution $(m,R)$ of the system \eqref{detailed eq:G-N system 1st eq}-\eqref{detailed eq:G-N system 2nd eq} belongs to $W^{2,p}_{\loc}(B^n)$ for any $p\in (1,\infty)$. 

The idea for the proof of Theorem \ref{thm:Lp regularity}  dates back to Sharp and Topping~\cite{Sharp-Topping-2013-TAMS} but with extra modifications from the recent works~\cite{Sharp-2014,Guo-Xiang-Zheng-2021-CV,Guo-Wang-Xiang-2023-CVPDE}. We will closely follow the presentation by Guo-Wang-Xiang~\cite{Guo-Wang-Xiang-2023-CVPDE} using a finite iteration method. The extra constant 1 appearing on the right-hand side of \eqref{eq:f R-m W^2,p estimate} comes from the fact that 
$$|\nabla m|-1\lesssim |S(\nabla m,R)|\lesssim |\nabla m|+1$$ 
and thus it cannot be removed. 

This paper is organized as follows. After the introduction, we collect all the necessary auxiliary results in Section \ref{sec:preliminaries}. In Section \ref{sec:Interior Regularity}, we prove Theorem \ref{thm:interior regularity} and in Section \ref{sec:L^p Regularity Theory}, we prove Theorem \ref{thm:Lp regularity}.

Our notations are standard. By $A\lesssim B$ we mean there exists a universal constant $C>0$ such that $A\le CB$.

\section{Preliminaries}\label{sec:preliminaries}
\subsection{Operators on matrices}\label{subsec:operators on matrices}
For a matrix $A=(a_{ij})_{N\times N}\in\R^{N\times N}$, we denote $A=(A_1|\cdots|A_N)$, where $A_i$ are the column vectors. The projection operator $\pi_n:\R^{N\times N}\to\R^{N\times n}$ (on the first $n$ columns) is defined by
\begin{equation*}\label{eq:def for Pi}
	\pi_n(A):=\pi(A_1|\cdots|A_N)=(A_1|\cdots|A_n)=\left(\begin{matrix}
		a_{11}&\cdots&a_{1n}\\
		\vdots&\ddots&\vdots\\
		a_{N1}&\cdots&a_{Nn}
	\end{matrix}\right)_{N\times n}.
\end{equation*}
The operator $\mathbb{P}:\R^{N\times N}\to\R^{N\times N}$ is defined by
\begin{equation*}\label{eq:def for P}
	\mathbb{P}(A):=\sqrt{\mu_1}\mathrm{dev}\mathrm{sym}(A)+\sqrt{\mu_2}\mathrm{skew}(A)+\sqrt{\kappa}\frac{\mathrm{tr}(A)}{N}\mathbbm{1}_N,
\end{equation*}
where $\mu_1,\mu_2$ and $\kappa$ are some positive constants,  and  we use
$$\mathrm{sym}(A)=\frac12(A+A^{\mathrm{T}}),\qquad\mathrm{skew}(A)=\frac12(A-A^{\mathrm{T}})$$ to represent the symmetric and skew-symmetric parts of $A$, respectively; the first term $\mathrm{dev}\mathrm{sym}(A)$ in the definition of $\mathbb{P}(A)$ is thus defined as $$\mathrm{dev}\mathrm{sym}(A)=\frac12(A+A^{\mathrm{T}})-\frac{\mathrm{tr}(A)}{N}\mathbbm{1}_N=\mathrm{sym}(A)-\frac{\mathrm{tr}(A)}{N}\mathbbm{1}_N$$
so as to denote the trace-free deviatoric part. As a result, there holds $$A=\mathrm{dev}\mathrm{sym}(A)+\mathrm{skew}(A)+\frac{\mathrm{tr}(A)}{N}\mathbbm{1}_N=\mathrm{sym}(A)+\mathrm{skew}(A),$$ 
which is an orthonormal decomposition for $A$. 

Finally, we define
\begin{equation}\label{eq:def for S}
	S(\nabla m,R):=\pi_n\left(2R\mathbb{P}^2\left(R^{\mathrm{T}}(\nabla m|0)-\diag(\mathbbm{1}_{n,N-n})\right)\right),
\end{equation}
where $$\diag(\mathbbm{1}_{n,N-n})=\left[\begin{matrix}
	\mathbbm{1}_n&O_{n\times(N-n)}\\
	O_{(N-n)\times n}&O_{(N-n)\times(N-n)}
\end{matrix}\right].$$ 
For the convenience of reading, we will use $(\mathbbm{1}_n|0)$ to represent $\diag(\mathbbm{1}_{n,N-n})$ in the following part of this article.

Finally we introduce an operation $\circ:\R^{N\times n}\times\R^{N\times n}\to\R^{N\times N}$ by
$$
A_{N\times n}\circ B_{N\times n}:=\frac12 AB^{\mathrm{T}}\in\R^{N\times N}.
$$ 
Then we have 
\begin{equation}\label{eq:def for skew m s}
	\mathrm{skew}\left(\nabla m\circ S(\nabla m,R)\right):=\frac14\left(\nabla m \left(S(\nabla m,R)\right)^{\mathrm{T}}-S(\nabla m,R)\nabla m^{\mathrm{T}}\right).
\end{equation}

\subsection{Morrey spaces and Riesz operators}

Let $1\le p<\infty$ and $0\le s\le n$. The Morrey space $M^{p,s}(U)$
consists of functions $u\in L^{p}(U)$ such that
\[
\|u\|_{M^{p,s}(U)}\equiv\sup_{x\in U,0<r<\diam(U)}r^{-s/p}\|u\|_{L^{p}(B_{r}(x)\cap U)}<\infty.
\]
The weak Morrey space $M_{\ast}^{p,s}(U)$ consists of functions
$f\in L_{\ast}^{p}(U)$ such that
\[
\|f\|_{M_{\ast}^{p,s}(U)}\equiv\sup_{x\in U,0<r<\diam(U)}r^{-s/p}\|f\|_{L_{\ast}^{p}(B_{r}(x)\cap U)}<\infty.
\]
The space $M_{1}^{p,s}(U)$ consists of functions in $M^{p,s}(U)$ whose weak gradient belongs to $M^{p,s}(U)$.

We need the following H\"older's inequalities in Morrey spaces; see~\cite{Guo-Wang-Xiang-2023-CVPDE}.

\begin{proposition}\label{prop: Holder for Morrey space} Let $1\le p_{1},p_{2}\le\infty$
	and $0\leq q_{1},q_{2}\le n$ be such that
	\begin{eqnarray*}
		\frac{1}{p}=\frac{1}{p_{1}}+\frac{1}{p_{2}}\le1 & \text{and} & q=\frac{p}{p_{1}}q_{1}+\frac{p}{p_{2}}q_{2}.
	\end{eqnarray*}
	Then, there hold
	\begin{equation}
		\|fg\|_{M^{p,q}(U)}\le\|f\|_{M^{p_{1},q_{1}}(U)}\|g\|_{M^{p_{2},q_{2}}(U)}.\label{eq: Morrey-Holder}
	\end{equation}
	and
	\begin{equation}
		\|fg\|_{M_{\ast}^{p,q}(U)}\le\|f\|_{M_{\ast}^{p_{1},q_{1}}(U)}\|g\|_{M_{\ast}^{p_{2},q_{2}}(U)}.\label{eq: Morrey-Holder-2}
	\end{equation}
\end{proposition}


Let $I_{\alpha}(x)=c_{\alpha,n}|x|^{\alpha-n}$, $0<\alpha<n$, be the standard Riesz potentials in $\R^{n}$. The following two propositions are well-known; see Theorem 3.1, Proposition 3.2 and Proposition 3.1 of Adams~\cite{Adams-1975}.
\begin{proposition}[\cite{Adams-1975}]\label{riesz potential M p,lambda}
	Let $0<\alpha<n$ and $0\le\lambda<n$. For $1\le p<(n-\lambda)/\alpha$, $\frac{1}{\tilde{p}}=\frac{1}{p}-\frac{\alpha}{n-\lambda}$, we have 
	\begin{enumerate}[(1)]
		\item for every $1<p<(n-\lambda)/\alpha$,$$
		I_\alpha:M^{p,\lambda}(\mathbb{R}^n)\to M^{\tilde{p},\lambda}(\mathbb{R}^n)
		$$is a bounded linear operator;
		
		\item for $p=1$, $$
		I_\alpha:M^{1,\lambda}(\R^n)\to M_*^{\tilde{p},\lambda}(\R^n)
		$$is also a bounded linear operator.
	\end{enumerate}
\end{proposition}

\begin{proposition}[\cite{Adams-1975}]\label{riesz potential M 1,n-beta cap Lp}
	Let $0<\alpha<\beta\le n$ and $1<p<\infty$. There exists a constant $C=C(\alpha,\beta,n,p)>0$ such that for $f\in M^{1,n-\beta}(\R^n)\cap L^p(\R^n)$, there holds$$
	\|I_\alpha f\|_{L^{\frac{p\beta}{\beta-\alpha}}(\mathbb{R}^n)}\leq C\|f\|_{M^{1,n-\beta}(\mathbb{R}^n)}^{\frac{\alpha}{\beta}}\|f\|_{L^p(\mathbb{R}^n)}^{1-\frac{\alpha}{\beta}}.$$
\end{proposition}

\subsection{Hardy-BMO inequality}
As an application of the Hardy-BMO duality and div-curl lemma (see~\cite{CLMS-1993,Evans-1991,Schikorra-2010}), we have the following Hardy-BMO inequality.

\begin{lemma}[Hardy-BMO inequality]\label{Hardy-BMO duality by L^p-L^q-Morrey}
	For any $p\in (1,\infty)$ and $\alpha\in (1,n)$, there exists a constant $C=C(n,p,\alpha)>0$ such that the following holds: 
	
	(1). For all balls $B_r(x_0)\subset \R^n$, functions $a\in M_1^{\alpha,n-\alpha}(B_{2r}(x_0))$, $\Gamma\in L^q(B_{r}(x_0),\R^n)$, $b\in W_0^{1,p}\cap L^\infty(B_r(x_0))$ with $\frac{1}{p}+\frac{1}{q}=1$ and $\Div(\Gamma)=0$ in the weak sense on $B_r(x_0)$, we have 
	$$
	\left|\int_{B_r(x_0)}\langle\nabla a,\Gamma\rangle bdx\right|\le C\left\|\Gamma\right\|_{L^q(B_r(x_0))}\|\nabla b\|_{L^p(B_r(x_0))}\|\nabla a\|_{M^{\alpha,n-\alpha}(B_{2r}(x_0))}.
	$$
	
	(2). For all balls $B_r(x_0)\subset \R^n$, functions $\varphi\in C_0^\infty(B_r(x_0))$, $\Gamma\in L^q(B_{r}(x_0),\R^n)$, $b\in M_{1}^{\alpha,n-\alpha}\cap L^\infty(B_{2r}(x_0))$ with $\frac{1}{p}+\frac{1}{q}=1$ and $\Div(\Gamma)=0$ in the weak sense on $B_r(x_0)$, we have
	$$
	\left|\int_{B_r(x_0)}\langle\nabla \varphi,\Gamma\rangle bdx\right|\le C\left\|\Gamma\right\|_{L^q(B_r(x_0))}\|\nabla \varphi\|_{L^p(B_r(x_0))}\|\nabla b\|_{M^{\alpha,n-\alpha}(B_{2r}(x_0))}.
	$$
\end{lemma}

We shall use the following technical inequality. 
\begin{lemma}\label{embed Hardy cap M^1,n-2 into H^-1}
	If $F\in M^{2,n-2}(\R^n,\R^n)$ is a divergence-free vector field, $g\in W^{1,2}(\R^n)$ and $\nabla g\in M^{2,n-2}(\R^n)$, then for each compact $K\subset\R^n$, we have
	\begin{equation*}\label{eq: embed Hardy cap M^1,n-2 into H^-1}
		\|F\cdot\nabla g\|^2_{H^{-1}(K)}\le C\|F\|_{L^2}\|F\|_{M^{2,n-2}}\|\nabla g\|_{L^2}\|\nabla g\|_{M^{2,n-2}}.
	\end{equation*}
\end{lemma}
\begin{proof}
	By Corollary~1.8 in~\cite{Sharp-2014}, we have 
	\[
	\|F\cdot\nabla g\|^2_{H^{-1}(K)}\le C 	\|F\cdot\nabla g\|_{h^{1}}	\|F\cdot\nabla g\|_{M^{1,n-2}}.
	\]
	The result follows then from the div-curl lemma and H\"older's inequality in Morrey spaces. 
\end{proof}


\subsection{Hodge decomposition}


We shall use the following well-known Hodge decomposition for $L^p$ integrable vector fields; see for instance~\cite[Chapter 10.5]{Iwaniec-Martin-2001-book}.
\begin{lemma}[\cite{Iwaniec-Martin-2001-book}]\label{Hodge decomposition}
	Let $p\in(1,\infty)$. Every vector field $V\in L^p(B_r(x_0),\R^n)$ with  $B_r(x_0)\subset\R^n$, can be uniquely decomposed as
	$$
	V=\nabla a+\nabla^{\perp} b+h,
	$$
	where $a\in W^{1,p}(B_r(x_0)), b\in W^{1,p}_0(B_r(x_0),\bigwedge^2\R^n)$ with $db=0$, and $h\in C^{\infty}(B_r(x_0),\R^n)$ is harmonic. Moreover, we have the estimate$$
	\|a\|_{W^{1,p}(B_r(x_0))}+\|b\|_{W^{1,p}(B_r(x_0))}+\|h\|_{L^p(B_r(x_0))}\leq C\|V\|_{L^p(B_r(x_0))}.
	$$	
	Here $\nabla^\perp b:=(\delta b)^\sharp\in L^p(B_r(x_0),\R^n)$, where $\delta$ is the formal conjugate operator of $d$ and $\sharp$ is the sharp operator from $\bigwedge^1\R^n$ to $\R^n$.
\end{lemma}	

We will also use the following ``nonlinear Hodge decomposition" for a connection matrix in a certain Morrey space, proved by Rivi\`ere and Struwe~\cite{RS-2008}.

\begin{lemma}[\cite{RS-2008}]\label{thm Riviere-Struwe 2008}
	There exists $\e(n,N)>0$ such that for every $\Omega\in M^{2,n-2}(B^n,\R^N\otimes so(N))$ with
	$$
	\|\Omega\|_{M^{2,n-2}(B^n)}^2<\e(n,N),
	$$
	there exist $P\in W^{1,2}(B^n,SO(N))$, $\xi\in W^{1,2}(B^n, so(N)\otimes \bigwedge^2\R^n)$ such that
	$$
	-P^{-1}\nabla P+P^{-1}\Omega P=\nabla^{\perp}\xi \quad\text{in }B^n
	$$
	and
	$$
	\|\nabla P\|_{M^{2,n-2}(B^n)}^2+\|\nabla \xi \|_{M^{2,n-2}(B^n)}^2\le C\|\Omega\|_{M^{2,n-2}(B^n)}^2\le C\e(n,N).
	$$
\end{lemma}

\section{Interior regularity}\label{sec:Interior Regularity}

In this section, we shall prove Theorem \ref{thm:interior regularity}. The key step towards it is the following H\"older continuity for weak solutions of \eqref{detailed eq:G-N system 1st eq}-\eqref{detailed eq:G-N system 2nd eq}.

\begin{theorem}\label{thm:local Holder continuity}
	There exists some $\epsilon=\epsilon(n,N)>0$ such that if $(m,R)\in W^{1,2}(B^n,\R^N\times SO(N))$ is a weak solution of \eqref{detailed eq:G-N system 1st eq}-\eqref{detailed eq:G-N system 2nd eq} with
	$$
	\|\nabla m\|_{M^{2,n-2}(B^n)}\le\e\quad \text{and}\quad \|\nabla R\|_{M^{2,n-2}(B^n)}\le\e,
	$$ 
	then there exists $\beta>0$ such that $m$ and $R$ are $C^{0,\beta}$-H\"older continuous on $B_{1/2}$.
\end{theorem}

\begin{proof}
	Fix any $x_0\in B^{n}$ and write $B_r$ for the ball $B_r(x_0)\subset B^n$, where $r$ is small enough such that $B_{2r}(x_0)\subset B^{n}$. Assume $\e\in(0,\e_0)$ for some $\e_0=\e_0(n,N)>0$ to be determined later. 
	
	Note that
	$$\quad\|\nabla m\|_{M^{2,n-2}(B^n)}\le\e \quad\text{and}\quad \|\Omega_R\|_{M^{2,n-2}(B^n)}\approx \|\nabla R\|_{M^{2,n-2}(B^n)}\le \e.$$
	According to Lemma \ref{thm Riviere-Struwe 2008}, there exist $P\in W^{1,2}(B^n,SO(N))$, $\xi\in W^{1,2}(B^n, so(N)\otimes \bigwedge^2\R^n)$ such that
	$$
	-P^{-1}\nabla P+P^{-1}\Omega_R P=\nabla^{\perp}\xi \quad\text{in }B^n
	$$
	and
	$$
	\|\nabla P\|_{M^{2,n-2}(B^n)}^2+\|\nabla \xi \|_{M^{2,n-2}(B^n)}^2\le C\|\Omega_R\|_{M^{2,n-2}(B^n)}^2\le C\e(n,N).
	$$
	Then direct computation shows
	\begin{equation}\label{eq: rewrite G-N equ}
		\begin{aligned}
			&\quad\Div(P^{-1}\nabla R) =  \nabla(P^{-1})\cdot\nabla R+P^{-1}\Delta R\\
			&= -P^{-1}(\nabla P)P^{-1}\cdot\nabla R+P^{-1}\Omega_{R}\cdot\nabla R+P^{-1}\mathrm{skew} \left(\nabla m\circ S(\nabla m,R)\right)R\\
			&= P^{-1}\mathrm{skew} \left(\nabla m\circ S(\nabla m,R)\right)R+(\nabla^\perp\xi)P^{-1}\nabla R.
		\end{aligned}
	\end{equation}
	
	By the Hodge decomposition and Lemma \ref{Hodge decomposition}, we may find functions $a\in W^{1,2}(B_r,\R^{N\times N})$, $b\in W_0^{1,2}(B_r,\R^{N\times N}\otimes\bigwedge^2\R^n)$ and a component-wise harmonic $h\in C^{\infty}(B_r,\R^{N\times N}\otimes\R^n)$ such that
	\begin{equation*}\label{eq: G-N eq2 Hodge decom}
		P^{-1}\nabla R=\nabla a+\nabla^\perp b+h\quad\text{in }B_r.
	\end{equation*}
	Define an operator $\Div^{\perp}:\R^n\to\bigwedge^2\R^n$ as $\Div^{\perp}V:=d(V^\flat)$, where $\flat$ is the flat operator from $\R^n$ to $\bigwedge^1\R^n$. Then for any $b\in\bigwedge^2\R^n$ with $db=0$ we have $$
	\Div^{\perp}\nabla^\perp b=d(\nabla^\perp b)^\flat=d\ast d\ast b=\Delta_{\R^n} b,
	$$$$
	\Div^{\perp}\nabla = d\circ\flat\circ\sharp\circ d=d\circ d=0
	$$and$$
	\Div\nabla^\perp = \Div\circ(\sharp\ast d\ast)=\ast d\ast\flat\sharp\ast d\ast=0.$$
	
	Using operators $\Div$ and $\Div^\perp$ to act on both sides of the equation, we obtain
	\begin{equation}\label{eq: a}
		\Delta a = \Div(P^{-1}\nabla R) = P^{-1}\mathrm{skew} \left(\nabla m\circ S(\nabla m,R)\right)R+(\nabla^\perp\xi)P^{-1}\nabla R
	\end{equation}
	and
	\begin{equation}\label{eq: b}
		\Delta b = d(P^{-1}dR) = dP^{-1}\wedge d(R-R_0),
	\end{equation}
	for any constant $R_0\in\R^{N\times N}$.
	
	To ease our notation, we set $\lambda=\frac{n+1}{n}$, whose H\"older conjugate index is $n+1$. In order to estimate $\|\nabla a\|_{L^{\lambda}(B_r)}$ and $\|\nabla^\perp b\|_{L^{\lambda}(B_r)}$, we define 
	$$T:=\left\{ \varphi\in C_0^\infty(B_r,\R^{N\times N}):\|\nabla \varphi\|_{L^{n+1}(B_r)}\le1\right\}.$$
	Then
	\begin{equation}\label{eq: a preliminary estimate}
		\begin{aligned}
			\|\nabla a\|_{L^{\lambda}(B_r)}\lesssim& \sup_{\varphi\in T}\int_{B_r}\left\langle\nabla a,\nabla \varphi\right\rangle dx\\
			\lesssim& \sup_{\varphi\in T}\int_{B_r}\left\langle P^{-1}\mathrm{skew} \left(\nabla m\circ S(\nabla m,R)\right)R+(\nabla^\perp\xi)P^{-1}\nabla R,- \varphi\right\rangle dx.
		\end{aligned}
	\end{equation}
	Note that $\frac{n-\lambda}{\lambda}=\frac{n-\frac{n+1}{n}}{\frac{n+1}{n}}=\frac{n^2-n-1}{n+1}$, $\frac{2-\lambda}{2\lambda}=\frac{n-1}{2n+2}$. Write $\xi=(\xi_{ij})$, $P^{-1}=(P^{-1}_{ij})$, $R=(R_{ij})$ and $\varphi=(\varphi_{ij})$ and observe that
	\[
	\left\langle(\nabla^\perp\xi)P^{-1}\nabla R, \varphi\right\rangle=\left\langle \nabla^\perp\xi_{ij},\nabla R_{kl} \right\rangle P_{jk}^{-1}\varphi_{il}.
	\]
	Since $\Div(\nabla^\perp \xi_{ij})=0$, we may apply Lemma \ref{Hardy-BMO duality by L^p-L^q-Morrey} (1) with $\Gamma=\nabla^\perp \xi_{ij}\in L^2$, $a=R_{kl}\in M_1^{\lambda,n-\lambda}$ and $b=P_{jk}^{-1}\varphi_{il}\in W^{1,2}_0\cap L^\infty$ to obtain 
	\begin{equation}\label{eq: a 1st term}
		\begin{aligned}
			&\int_{B_r}\left\langle(\nabla^\perp\xi)P^{-1}\nabla R, \varphi\right\rangle dx\\
			\lesssim& \|\nabla R\|_{M^{\lambda,n-\lambda}(B_{2r})}\|\nabla^\perp\xi\|_{L^2(B_r)} \|\nabla (P^{-1})\|_{L^2(B_r)}\|\varphi\|_{L^\infty(B_r)}\\
			&+\|\nabla R\|_{M^{\lambda,n-\lambda}(B_{2r})}\|\nabla^\perp\xi\|_{L^2(B_r)} \|\nabla\varphi\|_{L^{n+1}(B_r)}\|1\|_{L^{\frac{2n+2}{n-1}}(B_r)}\|P^{-1}\|_{L^\infty(B_r)}\\	
			\lesssim& \|\nabla R\|_{M^{\lambda,n-\lambda}(B_{2r})}\cdot r^{\frac{n-2}{2}}\e\cdot r^{\frac{n-2}{2}}\e\cdot\|\nabla \varphi\|_{L^{n+1}(B_r)}r^{1-\frac{n}{n+1}}\\
			&+ \|\nabla R\|_{M^{\lambda,n-\lambda}(B_{2r})}\cdot r^{\frac{n-2}{2}}\e\cdot r^{\frac{n(2-\lambda)}{2\lambda}}\|\nabla\varphi\|_{L^{n+1}(B_r)}\\
			\lesssim& \e r^{\frac{n-\lambda}{\lambda}}\|\nabla R\|_{M^{\lambda,n-\lambda}(B_{2r})}\|\nabla\varphi\|_{L^{n+1}(B_r)},
		\end{aligned}
	\end{equation}
	where in the second inequality we used the estimates
	\[
	\|\nabla^\perp \xi\|_{L^2(B_r)}\lesssim r^{\frac{n-2}{2}}\|\nabla \xi\|_{M^{2,n-2}(B_r)}\lesssim \e r^{\frac{n-2}{2}},
	\]
	\[
	\|\nabla (P^{-1})\|_{L^2(B_r)}\lesssim \|\nabla P\|_{L^2(B_r)}\lesssim r^{\frac{n-2}{2}}\|\nabla P\|_{M^{2,n-2}(B_r)}\lesssim \e r^{\frac{n-2}{2}}.
	\]
	
	Similarly, since $\Div(S(\nabla m,R))=0$, we may apply Lemma \ref{Hardy-BMO duality by L^p-L^q-Morrey} (1) with $\Gamma=S(\nabla m,R)_{ks}\in L^2$, $a=m_{js}\in M_1^{\lambda,n-\lambda}$ and $b=P_{ij}^{-1}R_{kl}\varphi_{il}\in W^{1,2}_0\cap L^\infty$ to obtain 
	\begin{equation}\label{eq: a 2nd term-1}
		\begin{aligned}
			&\quad \int_{B_r}\left\langle P^{-1}\mathrm{skew} \left(\nabla m\circ S(\nabla m,R)\right)R, \varphi\right\rangle dx\\
			&\lesssim \|\nabla m\|_{M^{\lambda,n-\lambda}(B_{2r})}\|S(\nabla m,R)\|_{L^2(B_r)}\|\nabla P\|_{L^2(B_r)}\|\varphi\|_{L^\infty(B_r)}\\
			&\quad+ \|\nabla m\|_{M^{\lambda,n-\lambda}(B_{2r})}\|S(\nabla m,R)\|_{L^2(B_r)}\|\nabla\varphi\|_{L^2(B_r)}\\
			&\quad+ \|\nabla m\|_{M^{\lambda,n-\lambda}(B_{2r})}\|S(\nabla m,R)\|_{L^2(B_r)}\|\nabla R\|_{L^2(B_r)}\|\varphi\|_{L^\infty(B_r)}\\
			&\lesssim r^{\frac{n(2-\lambda)}{2\lambda}}\|\nabla m\|_{M^{\lambda,n-\lambda}(B_{2r})}\|S(\nabla m,R)\|_{L^2(B_r)}\|\nabla\varphi\|_{L^{n+1}(B_r)}.
		\end{aligned}
	\end{equation}
	
	Note that the definition of $S(\nabla m,R)$ implies that there exists a constant $C$ such that
	\begin{equation}\label{eq: S estimate 0}
		C^{-1}\left( |\nabla m|-1\right)\le|S(\nabla m,R)|\le C\left( |\nabla m|+1\right).
	\end{equation}
	Then it follows
	\begin{equation}\label{eq: S estimate 1}
		\begin{aligned}
			\|S(\nabla m,R)\|_{L^2(B_r)}&\le C\left(\|\nabla m\|_{L^2(B_r)}+\|\mathbbm{1}\|_{L^2(B_r)}\right)\\
			&\le C\left(r^{\frac{n-2}{2}}\|\nabla m\|_{M^{2,n-2}(B_{2r})}+r^{\frac{n}{2}}\right)\\
			&\le Cr^{\frac{n-2}{2}}(r+\e),
		\end{aligned}
	\end{equation}
	and
	\begin{equation}\label{eq: S estimate 2}
		\begin{aligned}
			C^{-1}\|\nabla m\|_{L^{\lambda}(B_r)}-Cr^{\frac{n}{\lambda}}\le\|S(\nabla m,R)\|_{L^{\lambda}(B_r)}\le C\|\nabla m\|_{L^{\lambda}(B_r)}+Cr^{\frac{n}{\lambda}}.
		\end{aligned}
	\end{equation}
	Substituting \eqref{eq: S estimate 1} into \eqref{eq: a 2nd term-1} gives 
	\begin{equation}\label{eq: a 2nd term-2}
		\begin{aligned}
			&\quad\int_{B_r}\left\langle P^{-1}\mathrm{skew} \left(\nabla m\circ S(\nabla m,R)\right)R, \varphi\right\rangle dx\\
			&\le Cr^{\frac{n-\lambda}{\lambda}}(r+\e)\|\nabla m\|_{M^{\lambda,n-\lambda}(B_{2r})}\|\nabla\varphi\|_{L^{n+1}(B_r)}.
		\end{aligned}
	\end{equation}
	
	For convenience, we define $$
	\Phi(x_0,s):=\|\nabla m\|_{M^{\lambda,n-\lambda}(B_{s}(x_0))}+\|\nabla R\|_{M^{\lambda,n-\lambda}(B_{s}(x_0))}.$$
	Combining \eqref{eq: a preliminary estimate},  \eqref{eq: a 1st term} with \eqref{eq: a 2nd term-2}, we conclude 
	\begin{equation}\label{eq: a final estimate}
		\|\nabla a\|_{L^{\lambda}(B_r)}\le Cr^{\frac{n-\lambda}{\lambda}}(r+\e)\Phi(x_0,2r).
	\end{equation}

	Now we estimate the $L^{n+1}$-norm of $\nabla^\perp b$ using \eqref{eq: b}. Suppose $b=\sum_{1\le s<t\le n}b_{st}dx^s\wedge dx^t$, $dP^{-1}=\partial_sP^{-1} dx^s$, $dR=\partial_tR dx^t$. Then we have $$
	\Delta b_{st}= \partial_sP^{-1}\partial_tR-\partial_tP^{-1}\partial_sR=\nabla_{st}^\perp P^{-1}\cdot\nabla R,
	$$where $\nabla_{st}^\perp:=(0,\cdots,0,-\partial_t,0,\cdots,0,\partial_s,0,\cdots,0)$ for $1\le s<t\le n$. Note that $\Div\nabla_{st}^\perp P^{-1}=0$ for any $1\le s<t\le n$. Then we use Lemma \ref{Hardy-BMO duality by L^p-L^q-Morrey} (1) with $b=\varphi\in C_0^\infty$, $\Gamma=\nabla^{\perp}_{st}(P^{-1}_{ij})\in L^{\lambda}$ and $a=R_{jk}\in M_1^{\lambda,n-\lambda}\cap L^\infty$ as follows:
	$$
	\begin{aligned}
		\|\nabla b_{st}\|_{L^{\lambda}(B_r)}&\lesssim \sup_{\varphi\in T}\int_{B_r}\left\langle \nabla b_{st},\nabla \varphi \right\rangle dx\\
		&\lesssim \sup_{\varphi\in T}\int_{B_r}-\Delta b_{st}\varphi dx\\
		&\lesssim \sup_{\varphi\in T}\int_{B_r}-\nabla_{st}^\perp P^{-1}\cdot\nabla R\varphi dx\\
		&\lesssim \sup_{\varphi\in T} \|\nabla R\|_{M^{\lambda,n-\lambda}(B_{2r})}\|\nabla\varphi\|_{L^{n+1}(B_r)}\|\nabla P\|_{L^{\lambda}(B_r)}\\
		&\lesssim \e r^{\frac{n-\lambda}{\lambda}}\|\nabla R\|_{M^{\lambda,n-\lambda}(B_{2r})},
	\end{aligned}
	$$where in the last inequality we used the estimate  
	$$\|\nabla P\|_{L^{\lambda}(B_r)}\lesssim r^{\frac{n(2-\lambda)}{2\lambda}}\|\nabla P\|_{L^2(B_r)}\lesssim r^{\frac{n(2-\lambda)}{2\lambda}+\frac{n-2}{2}}\|\nabla P\|_{M^{2,n-2}(B_r)}\lesssim \e r^{\frac{n-\lambda}{\lambda}}.$$
	Thus \begin{equation}\label{eq: b final eatimate}
		\|\nabla^\perp b\|_{L^{\lambda}(B_r)}\lesssim\e r^{\frac{n-\lambda}{\lambda}}\|\nabla R\|_{M^{\lambda,n-\lambda}(B_{2r})}.
	\end{equation}

	For the harmonic term $h$, a standard estimate for harmonic functions implies that for any $0<\rho<r$, there holds
	\begin{equation}\label{eq: h final estimate}
		\int_{B_\rho}|h|^{\lambda}dx\le C\left(\frac{\rho}{r}\right)^n\int_{B_r}|h|^{\lambda}dx. 
	\end{equation}
	Combining \eqref{eq: a final estimate}, \eqref{eq: b final eatimate} with \eqref{eq: h final estimate}, we infer
	\begin{equation}\label{eq: R estimate}
		\begin{aligned}
			\|\nabla R\|_{L^{\lambda}(B_\rho)}&\lesssim  \|P^{-1}\nabla R\|_{L^{\lambda}(B_\rho)}\\
			&\lesssim  \|h\|_{L^{\lambda}(B_\rho)}+\|\nabla a\|_{L^{\lambda}(B_\rho)}+\|\nabla^\perp b\|_{L^{\lambda}(B_\rho)}\\
			&\lesssim  \left(\frac{\rho}{r}\right)^{\frac{n}{\lambda}}\|h\|_{L^{\lambda}(B_r)}+\|\nabla a\|_{L^{\lambda}(B_r)}+\|\nabla^\perp b\|_{L^{\lambda}(B_r)}\\
			&\lesssim  \left(\frac{\rho}{r}\right)^{\frac{n}{\lambda}}\|\nabla R\|_{L^{\lambda}(B_r)}+r^{\frac{n-\lambda}{\lambda}}(r+\e)\Phi(x_0,2r)\\
			&\quad+\e r^{\frac{n-\lambda}{\lambda}}\|\nabla R\|_{M^{\lambda,n-\lambda}(B_{2r})}\\	
			&\lesssim  \left(\frac{\rho}{r}\right)^{\frac{n}{\lambda}}\|\nabla R\|_{L^{\lambda}(B_r)}+(r+\e)r^{\frac{n-\lambda}{\lambda}}\Phi(x_0,2r).
		\end{aligned}
	\end{equation}
	
	It remains to estimate $\|\nabla m\|_{L^{\lambda}(B_r)}$. Applying the Hodge decomposition to the divergence-free matrix $S(\nabla m,R)\in\R^{N\times1}\otimes\R^n$, we obtain 
	\begin{equation}\label{eq: G-N eq1 Hodge decom}
		\pi_n\left(2R\mathbb{P}^2\left(R^{\mathrm{T}}(\nabla m|0)-(\mathbbm{1}_n|0)\right)\right)=\nabla^\perp\alpha+\chi,
	\end{equation} 
	where $\alpha\in W_0^{1,2}(B_r,\R^{N\times 1}\otimes\bigwedge^2\R^n)$ and $\chi\in C^\infty(B_r,\R^{N\times 1}\otimes\R^n)$ is (component-wise) harmonic. Due to \eqref{detailed eq:G-N system 1st eq}, we do not have terms of the form $\nabla\zeta$ on the right side. 
	
	For convenience, we define a linear map $\mathbb{P}_R$ as $\xi\mapsto2R\mathbb{P}(R^{\mathrm{T}}\xi)$. Using the operator $\Div^\perp$ to act on both sides of the equation, we obtain 
	\begin{equation}\label{eq: alpha eq}
		\begin{aligned}
			\Delta\alpha &= \Div^\perp\nabla^{\perp}\alpha\\
			&=\Div^\perp\left[ \pi_n(\mathbb{P}_R(\nabla m|0)-2R\mathbb{P}^{2}(\mathbbm{1}_{n}|0))\right]\\
			&=\pi_n(d\mathbb{P}_R\wedge(dm|0))-\Div^\perp\left[\pi_n(2R\mathbb{P}^{2}(\mathbbm{1}_{n}|0))\right].
		\end{aligned}
	\end{equation}
	Similar to the previous case, set $U:=\left\{ \psi\in C_0^\infty(B_r,\R^{N\times1}\otimes\bigwedge^2\R^n):\|\nabla^\perp  \psi\|_{L^{n+1}(B_r)}\le1\right\}$. It follows then from Lemma \ref{Hardy-BMO duality by L^p-L^q-Morrey} (2) and Poincar\'e's inequality that 
	{\begin{equation}\label{eq: alpha estimate}
			\begin{aligned}
				&\quad\|\nabla^\perp\alpha\|_{L^{\lambda}(B_r)}\lesssim \sup_{\psi\in U}\int_{B_r}\left\langle \nabla^\perp\alpha,\nabla^\perp\psi \right\rangle dx\\
				&= \sup_{\psi\in U}\int_{B_r}\left({\left\langle (d\mathbb{P}_R)m,\nabla^\perp\psi \right\rangle}
				-\left\langle \pi_n(2(R-R_{B_r})\mathbb{P}^{2}(\mathbbm{1}_{n}|0)), \nabla^\perp\psi \right\rangle\right)dx\\
				&\lesssim \sup_{\psi\in U}\|\nabla^\perp\psi\|_{L^{n+1}(B_r)}{\|d\mathbb{P}_R\|_{L^{\lambda}(B_r)}}\|\nabla m\|_{M^{\lambda,n-\lambda}(B_{2r})}\\
				&\quad + \sup_{\psi\in U}\|\nabla^\perp\psi\|_{L^{n+1}(B_r)}\|R-R_{B_r}\|_{L^{\lambda}(B_r)}\\
				&\lesssim \sup_{\psi\in U}\left[\|\nabla^\perp\psi\|_{L^{n+1}(B_r)}\left(\|\nabla R\|_{L^{\lambda}(B_r)}\|\nabla m\|_{M^{\lambda,n-\lambda}(B_{2r})}+ r\|\nabla R\|_{L^{\lambda}(B_r)}\right)\right]\\
				&\lesssim \left(\e r^{\frac{n-\lambda}{\lambda}}\|\nabla m\|_{M^{\lambda,n-\lambda}(B_{2r})}+ r^{\frac{n}{\lambda}}\|\nabla R\|_{M^{\lambda,n-\lambda}(B_{2r})}\right)\\
				&\lesssim (r+\e)r^{\frac{n-\lambda}{\lambda}}\Phi(x_0,2r),
			\end{aligned}
	\end{equation}}
	where in the third inequality we used the estimate  
	$$\|\nabla R\|_{L^{\lambda}(B_r)}\lesssim r^{\frac{n(2-\lambda)}{2\lambda}}\|\nabla R\|_{L^2(B_r)}\lesssim r^{\frac{n(2-\lambda)}{2\lambda}+\frac{n-2}{2}}\|\nabla R\|_{M^{2,n-2}(B_r)}\lesssim \e r^{\frac{n-\lambda}{\lambda}}.$$
	
	Returning to \eqref{eq: S estimate 2}, by \eqref{eq: alpha estimate} we have
	\begin{equation}\label{eq: m estimate}
		\begin{aligned}
			\|\nabla m\|_{L^{\lambda}(B_\rho)}&\lesssim \|\chi\|_{L^{\lambda}(B_\rho)}+\|\nabla^\perp\alpha\|_{L^{\lambda}(B_\rho)}+\rho^{\frac{n}{\lambda}} \\
			&\lesssim\left(\frac{\rho}{r}\right)^{\frac{n}{\lambda}}\|\chi\|_{L^{\lambda}(B_r)}+ \|\nabla^\perp\alpha\|_{L^{\lambda}(B_r)}+\rho^{\frac{n}{\lambda}} \\
			&\lesssim\left(\frac{\rho}{r}\right)^{\frac{n}{\lambda}}\|\nabla m\|_{L^{\lambda}(B_r)}+(\e+r)r^{\frac{n-\lambda}{\lambda}}\Phi(x_0,2r)+\rho^{\frac{n}{\lambda}}.
		\end{aligned}
	\end{equation}
	
	Combining \eqref{eq: R estimate} with \eqref{eq: m estimate}, we conclude 
	\begin{equation*}\label{eq: R-m La estimate 1}
		\begin{aligned}
			&\quad\rho^{-\frac{n-\lambda}{\lambda}}\left(\|\nabla m\|_{L^{\lambda}(B_\rho)}+\|\nabla R\|_{L^{\lambda}(B_\rho)}\right)\\
			&\lesssim \rho r^{-\frac{n}{\lambda}}\left(\|\nabla m\|_{L^{\lambda}(B_r)}+\|\nabla R\|_{L^{\lambda}(B_r)}\right)+ (\e+r)\left(\frac{r}{\rho}\right)^{\frac{n-\lambda}{\lambda}}\Phi(x_0,2r)+\rho\\
			&\lesssim \left(\frac{\rho}{r}+(\e+r)\left(\frac{r}{\rho}\right)^{\frac{n-\lambda}{\lambda}}\right)\Phi(x_0,2r)+\rho.
		\end{aligned}
	\end{equation*}
	
	We may assume $r\le\e\le \e_0$, where $\e_0$ will be determined soon. Adding $\rho$ on both sides leads to the estimate 
	\begin{equation}\label{eq: R-m La estimate 2}
		\begin{aligned}
			&\quad \rho^{-\frac{n-\lambda}{\lambda}}\left(\|\nabla m\|_{L^{\lambda}(B_\rho)}+\|\nabla R\|_{L^{\lambda}(B_\rho)}\right)+\rho\\
			&\le C_0\left(\frac{\rho}{r}+(\e+r)\left(\frac{r}{\rho}\right)^{\frac{n-\lambda}{\lambda}}\right)\left(\Phi(x_0,2r)+2r\right)
		\end{aligned}
	\end{equation}
	for suitable positive constant $C_0$ that depends only on $n,N$. 
	
	Now we fix $\rho=\theta r$, $\theta=\frac{1}{8C_0}$, and $\e_0=\frac12(8C_0)^{-\frac{n}{\lambda}}$. Then \eqref{eq: R-m La estimate 2} gives 
	\begin{equation}\label{eq: R-m La estimate 3}
		(\theta r)^{-\frac{n-\lambda}{\lambda}}\left(\|\nabla m\|_{L^{\lambda}(B_{\theta r})}+\|\nabla R\|_{L^{\lambda}(B_{\theta r})}\right)+\theta r\le \frac14\left(\Phi(x_0,2r)+2r\right),
	\end{equation}
	which holds for all $B_{\theta r}(x_0)$ and $B_{2r}(x_0)\subset B^n$. Clearly we can replace $B_{2r}(x_0)$ with any $B_s(y_0)\subset B^n$ containing $B_{2r}(x_0)$, provided that $s\le \e$.
	Thus \eqref{eq: R-m La estimate 3} implies 
	\begin{equation*}\label{eq: R-m La estimate 4}
		(\theta r)^{-\frac{n-\lambda}{\lambda}}\left(\|\nabla m\|_{L^{\lambda}(B_{\theta r}(x_0))}+\|\nabla R\|_{L^{\lambda}(B_{\theta r}(x_0))}\right)+\theta r\le \frac14\left(\Phi(y_0,s)+s\right),
	\end{equation*}
	which is valid for all $r,s,x_0,y_0$ such that $B_{2r}(x_0)\subset B_s(y_0)\subset  B^n$. Note that the family of balls $\{B_{\theta r}(x_0)\}$ forms an open cover of $B_{\frac{\theta s}{2}}(y_0)$. Thus we can take the supremum over all admissible $B_r(x_0)$ to find
	\begin{equation}\label{eq: R-m Morrey estimate1}
		\Phi(y_0,\frac{\theta s}{2})+\frac{\theta s}{2}\le \frac12\left(\Phi(y_0,s)+s\right).
	\end{equation}
	Setting $\Psi(y_0,r):=\Phi(y_0,r)+r$ and then iterating \eqref{eq: R-m Morrey estimate1}, we obtain 
	\begin{equation*}\label{eq: R-m iteration}
		\Psi(y_0,\left(\frac{\theta}{2}\right)^ks)\le 2^{-k}\Psi(y_0,s)\qquad \text{for all } k\in\N. 
	\end{equation*}
	For $r\approx\left(\frac{\theta}{2}\right)^ks$, we select $k\approx\frac{\log(r/s)}{\log(\theta/2)}$. Then $2^{-k}\approx (r/s)^{\frac{\log2}{\log(2/\theta)}}=:(r/s)^\beta$. This implies that for all $r\le s\le \e_0$, we have the estimate
	\begin{equation*}\label{eq: R-m estimate}
		\Psi(y_0,r)\le Cr^\beta s^{-\beta}\Psi(y_0,s).
	\end{equation*}
	We may choose $s=s_0>0$, depending only on $n,N$, such that for all $r\leq s_0/2$, there holds
	\begin{equation*}\label{eq: R-m final estimate}
		\Phi(y_0,r)\le \Psi(y_0,r)\le Cr^\beta s_0^{-\beta}\Psi(y_0,s_0),
	\end{equation*}
	which gives 
	\begin{equation*}\label{eq: R-m new Morrey space}
		\nabla m,\nabla R\in M_{\loc}^{\lambda,n-\lambda+\lambda\beta}(B^n).
	\end{equation*}
	Finally, Morrey's Dirichlet growth theorem (see for instance~\cite{Giaquinta-Book}) implies $(m,R)\in C_{\loc}^{0,\beta}$. The proof of Theorem \ref{thm:local Holder continuity} is thus complete.
\end{proof}

\begin{proof}[Proof of Theorem \ref{thm:interior regularity}]
	With Theorem \ref{thm:local Holder continuity} at hand, the proof follows directly from that of Gastel-Neff~\cite[Section 6.2]{Gastel-Neff 2022}.
\end{proof}


\section{$L^p$ regularity theory}\label{sec:L^p Regularity Theory}
In this section, we shall prove Theorem \ref{thm:Lp regularity}. The proof relies on the ideas from earlier works on similar problems in~\cite{Sharp-Topping-2013-TAMS,Sharp-2014,Guo-Xiang-Zheng-2021-CV,Guo-Wang-Xiang-2023-CVPDE}. In the first step, we show a quantitative H\"older continuity result for weak solutions of \eqref{detailed eq:G-N system 3rd eq}-\eqref{detailed eq:G-N system 4th eq}. 

\begin{proposition}\label{prop:2-n/p Holder}
	Under the same assumptions as Theorem \ref{thm:Lp regularity}, when $\frac{n}{2}<p<n$, we have $(m,R)\in C^{0,\gamma}(B_{1/2},\R^N\times SO(N))$, where $\gamma=2-n/p\in(0,1)$. Moreover, there exists some $C=C(n,N,p)>0$  such that 
	\begin{equation}\label{eq:f 2-n/p Holder quantitative estimate p<n}
		[m]_{C^{0,\gamma}(B_{1/2})}+[R]_{C^{0,\gamma}(B_{1/2})}\le C\left( \e+\|f\|_{L^p(B_1)}\right).
	\end{equation}
\end{proposition}

\begin{proof}
	Similar to the proof of Theorem \ref{thm:local Holder continuity}, we write $B_r$ for a fixed ball $B_r(x_0)\subset B_{1/2}$, where $r$ is small enough such that $B_{2r}(x_0)\subset B^n$. 
	
	Note that
	$$\quad\|\nabla m\|_{M^{2,n-2}(B^n)}\le\e \quad\text{and}\quad \|\Omega_R\|_{M^{2,n-2}(B^n)}\approx \|\nabla R\|_{M^{2,n-2}(B^n)}\le \e.$$
	According to Lemma \ref{thm Riviere-Struwe 2008}, there exist $P\in W^{1,2}(B^n,SO(N))$, {$\xi\in W^{1,2}(B^n, so(N)\otimes\bigwedge^2\R^n)$} such that
	$$
	-P^{-1}\nabla P+P^{-1}\Omega_R P=\nabla^{\perp}\xi \quad\text{in }B^n
	$$
	and
	$$
	\|\nabla P\|_{M^{2,n-2}(B^n)}^2+\|\nabla \xi \|_{M^{2,n-2}(B^n)}^2\le C\|\Omega_R\|_{M^{2,n-2}(B^n)}^2\le C\e(n,N).
	$$
	Straightforward computation gives 
	\begin{equation}\label{eq:f rewrite G-N equ}
		\begin{aligned}
			&\Div(P^{-1}\nabla R) =  \nabla(P^{-1})\cdot\nabla R+P^{-1}\Delta R\\
			=& -P^{-1}(\nabla P)P^{-1}\cdot\nabla R+P^{-1}\Omega_{R}\cdot\nabla R+P^{-1}\mathrm{skew} \left(\nabla m\circ S(\nabla m,R)\right)R+P^{-1}f\\
			=& P^{-1}\mathrm{skew} \left(\nabla m\circ S(\nabla m,R)\right)R+(\nabla^\perp\xi)P^{-1}\nabla R+P^{-1}f.
		\end{aligned}
	\end{equation}
	
	By the Hodge decomposition and Lemma \ref{Hodge decomposition}, there exist $a\in W^{1,2}(B_r,\R^{N\times N})$, {$b\in W_0^{1,2}(B_r,\R^{N\times N}\otimes\bigwedge^2\R^n)$} and a component-wise harmonic  $h\in C^{\infty}(B_r,\R^{N\times N}\otimes\R^n)$ such that
	\begin{equation}\label{eq:f G-N eq2 Hodge decom}
		P^{-1}\nabla R=\nabla a+\nabla^\perp b+h\quad\text{in }B_r.
	\end{equation}
	Then $a$ and $b$ satisfy 
	\begin{equation}\label{eq:f a}
		\Delta a = \Div(P^{-1}\nabla R) = P^{-1}\mathrm{skew} \left(\nabla m\circ S(\nabla m,R)\right)R+(\nabla^\perp\xi)P^{-1}\nabla R+P^{-1}f,
	\end{equation}
	and
	{\begin{equation}\label{eq:f b}
			\Delta b = \Div^\perp(P^{-1}\nabla R) = dP^{-1}\wedge dR = dP^{-1}\wedge d(R-R_0)
	\end{equation}}
	for any constant $R_0\in\R^{N\times N}$. 
	
	As in the proof of Theorem \ref{thm:local Holder continuity}, let $T:=\left\{ \varphi\in C_0^\infty(B_r,\R^{N\times N}):\|\nabla \varphi\|_{L^{n+1}(B_r)}\le1\right\}$ and set $\lambda=\frac{n+1}{n}$. Then, we have 
	\begin{equation*}\label{eq:f a preliminary estimate}
		\begin{aligned}
			&\quad \|\nabla a\|_{L^{\lambda}(B_r)}\lesssim \sup_{\varphi\in T}\int_{B_r}\left\langle\nabla f,\nabla \varphi\right\rangle dx\\
			&\lesssim \sup_{\varphi\in T}\int_{B_r}\left\langle P^{-1}\mathrm{skew} \left(\nabla m\circ S(\nabla m,R)\right)R+(\nabla^\perp\xi)P^{-1}\nabla R+P^{-1}f,-\varphi\right\rangle dx.
		\end{aligned}
	\end{equation*}
	
	Note that $n+\frac{1}{n+1}-\frac{n}{p}=\frac{n-\lambda}{\lambda}+2-\frac{n}{p}=\frac{n-\lambda}{\lambda}+\gamma$. Using H\"older's inequality and Sobolev's inequality, we deduce 
	\begin{equation}\label{eq:f P^-1f estimate}
		\begin{aligned}
			\int_{B_r}\left\langle P^{-1}f,\varphi\right\rangle dx&\lesssim \|P\|_{L^{\infty}(B_r)}\|f\|_{L^p(B_r)}\|1\|_{L^{\frac{p}{p-1}}(B_r)}\|\varphi\|_{L^{\infty}(B_r)}\\
			&\lesssim \|f\|_{L^p(B_r)}\cdot r^{n\left(1-\frac{1}{p}\right)}\cdot r^{\frac{1}{n+1}}\|\nabla\varphi\|_{L^{n+1}(B_r)}\\
			&\lesssim r^{\frac{n-\lambda}{\lambda}+\gamma}\|f\|_{L^p(B_r)}\|\nabla\varphi\|_{L^{n+1}(B_r)}.
		\end{aligned}
	\end{equation}
	
	As in the previous proof, we define $$
	\Phi(x_0,s):=\|\nabla m\|_{M^{\lambda,n-\lambda}(B_{s}(x_0))}+\|\nabla R\|_{M^{\lambda,n-\lambda}(B_{s}(x_0))}.$$
	Combining \eqref{eq:f P^-1f estimate} with \eqref{eq: a 1st term} and \eqref{eq: a 2nd term-2}, we infer
	\begin{equation}\label{eq:f a final estimate}
		\|\nabla a\|_{L^{\lambda}(B_r)}\le Cr^{\frac{n-\lambda}{\lambda}}\left((r+\e)\Phi(x_0,2r)+r^{\gamma}\|f\|_{L^p(B_r)}\right).
	\end{equation}
	
	Combining \eqref{eq:f a final estimate} with \eqref{eq: b final eatimate} and \eqref{eq: h final estimate}, we conclude
	\begin{equation*}\label{eq:f R final estimate}
		\|\nabla R\|_{L^{\lambda}(B_\rho)}
		\le C\left(\frac{\rho}{r}\right)^{\frac{n}{\lambda}}\|\nabla R\|_{L^{\lambda}(B_r)}+Cr^{\frac{n-\lambda}{\lambda}}\left((r+\e)\Phi(x_0,2r)+r^{\gamma}\|f\|_{L^p(B_r)}\right).
	\end{equation*}
	This together with \eqref{eq: m estimate} gives 
	\begin{equation}\label{eq:f R-m La estimate}
		\begin{aligned}
			&\rho^{-\frac{n-\lambda}{\lambda}}\left(\|\nabla m\|_{L^{\lambda}(B_\rho)}+|\nabla R\|_{L^{\lambda}(B_\rho)}\right)+\rho\\
			\le& C\rho r^{-\frac{n}{\lambda}}\left(\|\nabla m\|_{L^{\lambda}(B_r)}+|\nabla R\|_{L^{\lambda}(B_r)}\right)\\
			&+ C\left(\frac{r}{\rho}\right)^{\frac{n-\lambda}{\lambda}}\left((r+\e)\Phi(x_0,2r)+r^{\gamma}\|f\|_{L^p(B_r)}\right)+C\rho\\
			\le& C_0\left(\frac{\rho}{r}+(\e+r)\left(\frac{r}{\rho}\right)^{\frac{n-\lambda}{\lambda}}\right)\Phi(x_0,2r)+C_0\left(\frac{r}{\rho}\right)^{\frac{n-\lambda}{\lambda}}r^\gamma\|f\|_{L^p(B_1)}+C_0\rho
		\end{aligned}
	\end{equation}
	for some $C_0$ depending only on $n,N$. 
		
		Now select $\e$ small enough such that $(2\e)^{\frac{\lambda}{n}}=\frac12(2C_0)^{\frac{2}{\gamma-1}}$. Then let $\rho=\theta r$, $r\le\e$ and $\theta\in \left((2\e)^{\frac{\lambda}{n}},(2C_0)^{\frac{2}{\gamma-1}}\right)$. It follows
		$$
		C_0\left(\frac{\rho}{r}+(\e+r)\left(\frac{r}{\rho}\right)^{\frac{n-\lambda}{\lambda}}\right)\le C_0\theta\left(1+2\e\theta^{-\frac{n}{\lambda}}\right)\le 2C_0\theta\le \theta^{\frac{\gamma+1}{2}}.
		$$
		This together with \eqref{eq:f R-m La estimate} gives
		\begin{equation*}\label{eq:f R-m Morrey estimate 1}
			\Phi(x_0,\theta r)+\theta r\le \theta^{\frac{\gamma+1}{2}}\left(\Phi(x_0,2r)+2r\right)+C_1r^\gamma\|f\|_{L^p(B_1)},
		\end{equation*}
		where $C_1$ is a constant depending only on $n,N$. By a standard iteration argument,  we eventually obtain 
		\begin{equation*}
			\begin{aligned}
				\Phi(x_0,\theta r)+\theta r\le C\left(\theta^\gamma\Phi(x_0,2r)+2\theta^\gamma r+r^\gamma\|f\|_{L^p(B_1)}\right).
			\end{aligned}
		\end{equation*}
		This implies 
		\begin{equation*}
			\begin{aligned}
				\Phi(x_0,r)\le Cr^\gamma\left(\e+\|f\|_{L^p(B_1)}\right),
			\end{aligned}
		\end{equation*}
		from which we conclude that $\nabla m,\nabla R\in M_{\loc}^{\lambda,n-\lambda+\gamma\lambda}$.  By Morrey's Dirichlet growth theorem, this further implies that $m,R\in C^{0,\gamma}_{\loc}$ together with the desired estimate \eqref{eq:f 2-n/p Holder quantitative estimate p<n}.
	\end{proof}
	
	Next, we prove an improved Morrey regularity estimate, in the spirit of~\cite[Lemma 7.3]{Sharp-Topping-2013-TAMS} (or~\cite[Proposition 2.1]{Sharp-2014}).
	\begin{proposition}\label{prop:Morrey 2,n-2+2gamma}
		Under the same assumption as in Theorem \ref{thm:Lp regularity}, when $\frac{n}{2}<p<n$, we have $(m,R)\in M^{2,n-2+2\gamma}_{\loc}(B^n,\R^N\times SO(N))$, where $\gamma=2-n/p\in(0,1)$. Moreover, there exists some $C=C(n,N,p)>0$ such that
		\begin{equation}\label{eq:f R-m M^2,n-2+2gamma estimate}
			\|\nabla m\|_{M^{2,n-2+2\gamma}(B_{1/2})}+\|\nabla R\|_{M^{2,n-2+2\gamma}(B_{1/2})}\le C\left(\|f\|_{L^p(B_1)}+\e\right).
		\end{equation}
	\end{proposition}

	\begin{proof}
		We shall use \eqref{eq:f a}, \eqref{eq:f b} and \eqref{eq: alpha eq} to estimate $\|\nabla m\|_{M^{2,n-2+2\gamma}}$ and $\|\nabla R\|_{M^{2,n-2+2\gamma}}$.  First of all, we extend $m,R,P,\xi,f$ to $\R^n$ with compact support in a norm-bounded way. 
		\smallskip
		
		\noindent\textbf{Step 1.} Estimate $\|\nabla R\|_{L^2(B_{\rho})}$.\\
		
		Rewriting \eqref{eq:f a} as before, we have
		\begin{equation}\label{eq:f a rewrite}
			\begin{aligned}
				\Delta a &= P^{-1}\mathrm{skew} \left(\nabla m\circ S(\nabla m,R)\right)R+(\nabla^\perp\xi)\nabla\left( P^{-1} (R-R_{x_0,r})\right)\\
				&\quad-(\nabla^\perp\xi)\nabla\left( P^{-1} \right)(R-R_{x_0,r})+P^{-1}f,
			\end{aligned}
		\end{equation}
		where $R_{x_0,r}=\frac{1}{|B_r(x_0)|}\int_{B_r(x_0)}Rdx$.
		Now we calculate the $H^{-1}$-norm of the right-hand side of \eqref{eq:f a rewrite}. 
		
		We abbreviate $S(\nabla m,R)$ as $S$. For the first item on the  right-hand side of \eqref{eq:f a rewrite}, we have $$
		P^{-1}\mathrm{skew} \left(\nabla m\circ S\right)R = \frac14P^{-1}(\nabla m) S^{\mathrm{T}}R-\frac14P^{-1}S(\nabla m^{\mathrm{T}})R
		$$and$$
		\begin{aligned}
			P^{-1}(\nabla m) S^{\mathrm{T}}R&=\nabla\left(P^{-1}(m-m_{x_0,r})S^{\mathrm{T}}R\right)\\
			&\quad-(\nabla P^{-1})(m-m_{x_0,r})S^{\mathrm{T}}R-P^{-1}(m-m_{x_0,r})S^{\mathrm{T}}(\nabla R).
		\end{aligned}
		$$
		Thus, according to \eqref{eq: S estimate 0}, we may estimate the first term above as follows:
		\begin{equation}\label{eq:f a r.h.s estimate 1.1}
			\begin{aligned}	
				&\quad\|\nabla\left(P^{-1}(m-m_{x_0,r})S^{\mathrm{T}}R\right)\|_{H^{-1}(B_r)}\lesssim \|P^{-1}(m-m_{x_0,r})S^{\mathrm{T}}R\|_{L^2(B_r)}\\
				&\lesssim \|P^{-1}\|_{L^\infty(B_r)}\|(m-m_{x_0,r})\|_{L^\infty(B_r)}\|S\|_{L^2(B_r)}\|R\|_{L^\infty(B_r)}\\
				&\lesssim r^\gamma[m]_{C^{0,\gamma}(B_r)}\|S\|_{L^2(B_r)}
				\lesssim r^{\frac{n-2}{2}+\gamma}[m]_{C^{0,\gamma}(B_r)}\left(r+\|\nabla m\|_{M^{2,n-2}(B_r)}\right).
			\end{aligned}
		\end{equation}
		For the middle term, we may apply Lemma \ref{embed Hardy cap M^1,n-2 into H^-1} to obtaining
		\begin{equation}\label{eq:f a r.h.s estimate 1.2}
			\begin{aligned}
				&\quad\|(\nabla P^{-1})(m-m_{x_0,r})S^{\mathrm{T}}R\|_{H^{-1}(B_r)}\\
				&\lesssim \|R\|_{L^{\infty}(B_r)}\|(m-m_{x_0,r})\|_{L^{\infty}(B_r)}\|(\nabla P^{-1})S^{\mathrm{T}}\|_{H^{-1}(B_r)}\\
				&\lesssim r^\gamma[m]_{C^{0,\gamma}(B_r)}\|\nabla P^{-1}\|_{L^2(B_r)}^{\frac12}\|\nabla P^{-1}\|_{M^{2,n-2}(B_r)}^{\frac12}\|\nabla S\|_{L^2(B_r)}^{\frac12}\|\nabla S\|_{M^{2,n-2}(B_r)}^{\frac12}\\
				&\lesssim r^{\frac{n-2}{2}+\gamma}[m]_{C^{0,\gamma}(B_r)}\|\nabla R\|_{M^{2,n-2}(B_r)}\left(r+\|\nabla m\|_{M^{2,n-2}(B_r)}\right)
			\end{aligned}
		\end{equation}
		and
		\begin{equation}\label{eq:f a r.h.s estimate 1.3}
			\begin{aligned}
				&\quad\|P^{-1}(m-m_{x_0,r})S^{\mathrm{T}}(\nabla R)\|_{H^{-1}(B_r)}\\
				&\lesssim \|P^{-1}\|_{L^{\infty}(B_r)}\|(m-m_{x_0,r})\|_{L^{\infty}(B_r)}\|S^{\mathrm{T}}(\nabla R)\|_{H^{-1}(B_r)}\\
				&\lesssim r^\gamma[m]_{C^{0,\gamma}(B_r)}\|\nabla S\|_{L^2(B_r)}^{\frac12}\|\nabla S\|_{M^{2,n-2}(B_r)}^{\frac12}\|\nabla R\|_{L^2(B_r)}^{\frac12}\|\nabla R\|_{M^{2,n-2}(B_r)}^{\frac12}\\
				&\lesssim r^{\frac{n-2}{2}+\gamma}[m]_{C^{0,\gamma}(B_r)}\|\nabla R\|_{M^{2,n-2}(B_r)}\left(r+\|\nabla m\|_{M^{2,n-2}(B_r)}\right).
			\end{aligned}
		\end{equation}
		Combining \eqref{eq:f a r.h.s estimate 1.1}, \eqref{eq:f a r.h.s estimate 1.2} with \eqref{eq:f a r.h.s estimate 1.3}, we conclude
		\begin{equation}\label{eq:f a r.h.s estimate 1 final}
			\begin{aligned}
				&\quad\|P^{-1}\mathrm{skew} \left(\nabla m\circ S\right)R\|_{H^{-1}(B_r)}\\
				&\lesssim r^{\frac{n-2}{2}+\gamma}[m]_{C^{0,\gamma}(B_r)}\left(1+\|\nabla R\|_{M^{2,n-2}(B_r)}\right)\left(r+\|\nabla m\|_{M^{2,n-2}(B_r)}\right)\\
				&\lesssim r^{\frac{n-2}{2}+\gamma}(r+\e)[m]_{C^{0,\gamma}(B_r)}.
			\end{aligned}
		\end{equation}

		For the second item on the  right-hand side of \eqref{eq:f a rewrite}, we introduce 
		$$T:=\left\{ \varphi\in C_0^\infty(B_r,\R^{N\times N}):\|\nabla \varphi\|_{L^2(B_r)}\le1\right\}.$$ 
		Then we have
		\begin{equation*}\label{eq:f a r.h.s estimate 2.1}
			\begin{aligned}
				&\quad\left|\int_{B_r}(\nabla^\perp\xi)\nabla\left( P^{-1} (R-R_{x_0,r})\right)\varphi\right|=\left|-\int_{B_r}(\nabla^\perp\xi)P^{-1} (R-R_{x_0,r})\nabla\varphi\right|\\
				&\lesssim \|P^{-1}\|_{L^\infty(B_r)}\|\nabla^\perp
				\xi\|_{L^2(B_r)}\|\nabla\varphi\|_{L^2(B_r)}\|R-R_{x_0,r}\|_{L^\infty(B_r)}\\
				&\lesssim r^{\frac{n-2}{2}}\|\nabla^\perp
				\xi\|_{M^{2,n-2}(B_r)}\cdot r^\gamma[R]_{C^{0,\gamma}(B_r)}\cdot\|\nabla\varphi\|_{L^2(B_r)},
			\end{aligned}
		\end{equation*}
		which implies
		\begin{equation}\label{eq:f a r.h.s estimate 2 final}
			\|(\nabla^\perp\xi)\nabla\left( P^{-1} (R-R_{x_0,r})\right)\|_{H^{-1}(B_r)}\le Cr^{\frac{n-2}{2}+\gamma}\|\nabla R\|_{M^{2,n-2}(B_r)}[R]_{C^{0,\gamma}(B_r)}.
		\end{equation}
		
		For the third term on the right-hand side of \eqref{eq:f a rewrite}, we have 
		\begin{equation}\label{eq:f a r.h.s estimate 3 final}
			\begin{aligned}
				&\quad\|(\nabla^\perp\xi)\nabla\left( P^{-1} \right)(R-R_{x_0,r})\|_{H^{-1}(B_r)}\\
				&\lesssim \|(R-R_{x_0,r})\|_{L^{\infty}(B_r)}\|(\nabla^\perp\xi)\nabla\left( P^{-1} \right)\|_{H^{-1}(B_r)}\\
				&\lesssim r^\gamma[R]_{C^{0,\gamma}(B_r)}\|\nabla^\perp\xi\|_{L^2(B_r)}^{\frac12}\|\nabla^\perp\xi\|_{M^{2,n-2}(B_r)}^{\frac12}\|\nabla P\|_{L^2(B_r)}^{\frac12}\|\nabla P\|_{M^{2,n-2}(B_r)}^{\frac12}\\
				&\lesssim r^{\frac{n-2}{2}+\gamma}[R]_{C^{0,\gamma}(B_r)}\|\nabla^\perp\xi\|_{M^{2,n-2}(B_r)}\|\nabla P\|_{M^{2,n-2}(B_r)}\\
				&\lesssim r^{\frac{n-2}{2}+\gamma}[R]_{C^{0,\gamma}(B_r)}\|\nabla R\|_{M^{2,n-2}(B_r)}^2.
			\end{aligned}
		\end{equation}
		
		For the last term on the  right-hand side of \eqref{eq:f a rewrite}, we simply apply the embedding $L^{\frac{2n}{n+2}}\hookrightarrow H^{-1}$ to deriving 
		\begin{equation}\label{eq:f a r.h.s estimate 4 final}
			\|P^{-1}f\|_{H^{-1}(B_r)}\le C\|f\|_{L^{\frac{2n}{n+2}}(B_r)}\le Cr^{\frac{n}{2}+1-\frac{n}{p}}\|f\|_{L^p(B_r)}
			\le Cr^{\frac{n-2}{2}+\gamma}\|f\|_{L^p(B_r)}.
		\end{equation}

		Combining \eqref{eq:f a r.h.s estimate 1 final}, \eqref{eq:f a r.h.s estimate 2 final}, \eqref{eq:f a r.h.s estimate 3 final} with \eqref{eq:f a r.h.s estimate 4 final}, we conclude
		\begin{equation}\label{eq:f a L2 final estimate}
			\begin{aligned}
				\|\nabla a\|_{L^2(B_r)}\le Cr^{\frac{n-2}{2}+\gamma}\left[\|f\|_{L^p(B_r)}+(r+\e)\left([m]_{C^{0,\gamma}(B_r)}+[R]_{C^{0,\gamma}(B_r)}\right)\right],
			\end{aligned}
		\end{equation}
		where we used the assumption$$
		\|\Omega_R\|_{M^{2,n-2}(B_r)}\le C\|\nabla R\|_{M^{2,n-2}(B_r)}\le \e\quad\text{and} \quad\|\nabla m\|_{M^{2,n-2}(B_r)}\le\e.
		$$
		
		Now we estimate $\nabla b$ via \eqref{eq:f b}. {For any $2$-form $\varphi$ whose component belongs to $T$}, we have
		{
			\begin{equation*}\label{eq:f b L2 estimate 1}
				\begin{aligned}
					&\quad\left|\int_{B_r}\left\langle d\left((d P^{-1})(R-R_{x_0,r})\right),\varphi\right\rangle \right|= 	\left|\int_{B_r}(\left\langle d P^{-1})(R-R_{x_0,r}),\nabla^\perp\varphi\right\rangle \right|\\
					&\lesssim \|R-R_{x_0,r}\|_{L^\infty(B_r)}\|d P^{-1}\|_{L^2(B_r)}\|\nabla^\perp\varphi\|_{L^2(B_r)}\\
					&\lesssim r^{\frac{n-2}{2}+\gamma}[R]_{C^{0,\gamma}(B_r)}\|\Omega_R\|_{M^{2,n-2}(B_r)}\|\varphi\|_{W^{1,2}(B_r)}.
				\end{aligned}
		\end{equation*}}
		This gives 
		\begin{equation}\label{eq:f b L2 final estimate}
			\|\nabla b\|_{L^2(B_r)}\le C\e r^{\frac{n-2}{2}+\gamma}[R]_{C^{0,\gamma}(B_r)}.
		\end{equation}
		
		Combining \eqref{eq:f a L2 final estimate} and \eqref{eq:f b L2 final estimate} with the standard estimate for $h$, we conclude
		\begin{equation}\label{eq:f R L2 final estimate}
			\begin{aligned}
				&\quad\|\nabla R\|_{L^2(B_\rho)}\le \|h\|_{L^2{B_\rho}}+\|\nabla a\|_{L^2(B_\rho)}+\|\nabla^\perp b\|_{L^2(B_\rho)}\\
				&\lesssim \left(\frac{\rho}{r}\right)^{\frac{n}{2}}\|h\|_{L^2(B_r)}+\|\nabla a\|_{L^2(B_\rho)}+\|\nabla^\perp b\|_{L^2(B_\rho)}\\
				&\lesssim \left(\frac{\rho}{r}\right)^{\frac{n}{2}}\|\nabla R\|_{L^2(B_r)}+C\rho^{\frac{n-2}{2}+\gamma}\left[\|f\|_{L^p(B_\rho)}+(\rho+\e)\left([m]_{C^{0,\gamma}(B_\rho)}+[R]_{C^{0,\gamma}(B_\rho)}\right)\right].
			\end{aligned}
		\end{equation}
		\smallskip
		
		\noindent\textbf{Step 2.} Estimate $\|\nabla m\|_{L^2(B_{\rho})}$.
		\smallskip

		Applying the Hodge decomposition to the divergence-free matrix $S(\nabla m,R)$, we obtain \eqref{eq: G-N eq1 Hodge decom} as in the proof of Theorem \ref{thm:local Holder continuity}.
		
		To estimate $\nabla^\perp\alpha$, we shall use the following equation for $\alpha$:
		{\begin{equation}\label{eq:f alpha eq 3}
				\Delta\alpha =\pi_n(d\mathbb{P}_R\wedge(d(m-m_{x_0,r})|0))-\Div^\perp\left[\pi_n(2R\mathbb{P}^{2}(\mathbbm{1}_{n}|0))\right].
			\end{equation}
			Let $U:=\left\{ \psi\in C_0^\infty(B_r,\R^{N\times n}\otimes\bigwedge^2\R^n):\|\nabla^\perp \psi\|_{L^2(B_r)}\le1\right\}$.} Then we have
		\begin{equation*}\label{eq:f alpha estimate}
			\begin{aligned}
				&\quad\|\nabla^\perp\alpha\|_{L^2(B_r)}\lesssim \sup_{\psi\in U}\int_{B_r}\left\langle \nabla^\perp\alpha,\nabla^\perp\psi \right\rangle dx\\
				&\lesssim \sup_{\psi\in U}\int_{B_r}\left({\left\langle (d\mathbb{P}_R)(m-m_{x_0,r}),\nabla^\perp\psi \right\rangle} -\left\langle \pi_n(2(R-R_{x_0,r})\mathbb{P}^{2}(\mathbbm{1}_{n}|0)), \nabla^\perp\psi \right\rangle\right)dx\\
				&\lesssim \sup_{\psi\in U}\|\nabla^\perp\psi\|_{L^2(B_r)}\|d\mathbb{P}_R\|_{L^2(B_r)}\|m-m_{x_0,r}\|_{L^\infty(B_r)}\\
				&\quad+ \sup_{\psi\in U}\|\nabla^\perp\psi\|_{L^2(B_r)}\|R-R_{x_0,r}\|_{L^\infty(B_r)}\|1\|_{L^2(B_r)}\\
				&\lesssim \|\nabla R\|_{L^2(B_r)}\cdot r^\gamma[m]_{C^{0,\gamma}(B_r)}+ r^\gamma[R]_{C^{0,\gamma}(B_r)}\cdot r^{\frac{n}{2}}\\
				&\lesssim  r^{\frac{n-2}{2}+\gamma}(r+\e)\left([m]_{C^{0,\gamma}(B_r)}+[R]_{C^{0,\gamma}(B_r)}\right).
			\end{aligned}
		\end{equation*}
		
		Finally, we may can estimate $\nabla m$ as follows:
		\begin{equation}\label{eq:f m L2 final estimate}
			\begin{aligned}
				&\quad \|\nabla m\|_{L^2(B_\rho)}\lesssim \left(\|\chi\|_{B_\rho}+\|\nabla^\perp \alpha\|_{L^2(B_\rho)}+\rho^{\frac{n}{2}}\right)\\
				&\lesssim \left(\frac{\rho}{r}\right)^{\frac{n}{2}}\|\chi\|_{L^2(B_r)}+ \rho^{\frac{n-2}{2}+\gamma}(\rho+\e)\left([m]_{C^{0,\gamma}(B_\rho)}+[R]_{C^{0,\gamma}(B_\rho)}\right)+\rho^{\frac{n}{2}}\\
				&\lesssim
				\left(\frac{\rho}{r}\right)^{\frac{n}{2}}\|\nabla m\|_{L^2(B_r)}+ \rho^{\frac{n-2}{2}+\gamma}(\rho+\e)\left([m]_{C^{0,\gamma}(B_\rho)}+[R]_{C^{0,\gamma}(B_\rho)}\right)+\rho^{\frac{n}{2}}.
			\end{aligned}
		\end{equation}
		\smallskip
		
		\noindent\textbf{Step 3.} Iteration. 
		\smallskip 
		
		Define $\Psi(x_0,r)=\|\nabla m\|_{L^2(B_r)}+\|\nabla R\|_{L^2(B_r)}$. Combining \eqref{eq:f R L2 final estimate} with \eqref{eq:f m L2 final estimate}, we have
		\begin{equation*}\label{eq:f R-m L2 estimate 1}
			\Psi(x_0,\rho)\le C_0\left(\frac{\rho}{r}\right)^{\frac{n}{2}}\Psi(x_0,r)+C_0\rho^{\frac{n-2}{2}+\gamma}\Gamma,
		\end{equation*}
		where $$
		\begin{aligned}
			\Gamma&=(\rho+\e)\left([m]_{C^{0,\gamma}(B_\rho)}+[R]_{C^{0,\gamma}(B_\rho)}\right)+\|f\|_{L^p(B_\rho)}+\rho^{1-\gamma}\\
			&\leq C\left(\|f\|_{L^p(B_1)}+\e\right)=:\Gamma_0.
		\end{aligned}
		$$
		A standard iteration argument gives
		\begin{equation*}\label{eq:f R-m L2 estimate 3}
			\Psi(x_0,r)\le r^{\frac{n-2}{2}+\gamma}C\left(\|f\|_{L^p(B_1)}+\e\right),
		\end{equation*}
		from which we obtain \eqref{eq:f R-m M^2,n-2+2gamma estimate}. The proof of Proposition \ref{prop:Morrey 2,n-2+2gamma} is thus complete. 
	\end{proof}
	
	\smallskip
	\begin{proof}[Proof of Theorem \ref{thm:Lp regularity}]
		
		We shall consider two cases separately. 
		\smallskip 
		
		\textbf{Case 1.} $p\in(\frac{n}{2},n)$. 
		\smallskip 
		
		By Proposition \ref{prop:Morrey 2,n-2+2gamma} and H\"older's inequality  (see Proposition \ref{prop: Holder for Morrey space}), we have 
		$$
		\Omega_R\cdot\nabla R\in M_{\loc}^{1,n-2+2\gamma}\hookrightarrow M_{\loc}^{1,n-2+\gamma},
		$$
		$$
		\mathrm{skew} \left(\nabla m\circ S(\nabla m,R)\right)R \in M_{\loc}^{1,n-2+2\gamma}\hookrightarrow M_{\loc}^{1,n-2+\gamma}.
		$$
		Extend $m,R$ and $f$ from $B^n$ into $\R^n$ with compact support in a norm-bounded way. 
		
		Let $I_{\alpha}=c|x|^{\alpha-n}$ be  the standard Riesz potential. Set 
		$$R_1=I_2\left( \Omega_{R}\cdot\nabla R+\mathrm{skew} \left(\nabla m\circ S(\nabla m,R)\right)R\right)\quad\text{and}\quad R_2=I_2\left(f\right)$$ so that $R_3 = R-R_1-R_2$ is harmonic. 
		
		Note that $\frac{2-\gamma}{1-\gamma}>2$ and $\zeta=\frac12\left( \frac{2-\gamma}{1-\gamma}\right)>1$. Then Proposition \ref{riesz potential M p,lambda}, together with Propositions \ref{prop: Holder for Morrey space} and \ref{prop:Morrey 2,n-2+2gamma}, implies that for any $x_0,r$ such that $B_{2r}(x_0)\subset B_{1/2}$, we have
		\begin{equation}\label{eq:f R1 L^2zeta estimate 1}
			\begin{aligned}
				\|\nabla R_1\|_{M^{2\zeta,n-2+\gamma}_*(B_r)}&\lesssim \|\Omega_{R}\cdot\nabla R+\mathrm{skew} \left(\nabla m\circ S(\nabla m,R)\right)R)\|_{M^{1,n-2+\gamma}(B_{r})}\\
				&\lesssim \|\nabla R\|_{M^{2,n-2}(B_{r})}\|\nabla R\|_{M^{2,n-2+2\gamma}(B_{r})}\\
				&\quad+ \|\nabla m\|_{M^{2,n-2}(B_{r})}\|S(\nabla m,R)\|_{M^{2,n-2+2\gamma}(B_{r})}\\
				&\stackrel{\eqref{eq: S estimate 0}}{\lesssim} \e\|\nabla R\|_{M^{2,n-2+2\gamma}(B_{r})}+ \e\|\nabla m\|_{M^{2,n-2+2\gamma}(B_{r})}\\
				&\quad+\|\nabla m\|_{M^{2,n-2}(B_{r})}\|1\|_{M^{2,n-2+2\gamma}(B_{r})}\\
				&\lesssim \e\left( \|\nabla R\|_{M^{2,n-2+2\gamma}(B_{r})}+ \|\nabla m\|_{M^{2,n-2+2\gamma}(B_{r})}+1\right)\\
				&\lesssim \e\left( \|f\|_{L^p(B_1)}+1\right),
			\end{aligned}
		\end{equation}
		where in the fourth inequality we used the estimate 
		$$
		\|\nabla m\|_{M^{2,n-2}(B_{r})}\|1\|_{M^{2,n-2+2\gamma}(B_{r})}\le  \e r^{1-\gamma}\le \e.
		$$

		By the standard elliptic regularity theory, we have $R_2\in W^{2,p}(B_r)$ with 
		\begin{equation*}\label{eq:f R2 L^p* estimate 1}
			\|\nabla R_2\|_{L^{\frac{np}{n-p}}(B_r)}\le \| R_2\|_{W^{2,p}(B_r)}\le C\|f\|_{L^p(B_{1})}.
		\end{equation*}
		
		Applying the Hodge decomposition to the divergence-free matrix $S(\nabla m,R)$, {there exist $\alpha\in W_0^{1,2}(B_r,\R^{N\times 1}\otimes\bigwedge^2\R^n)$ and a harmonic $\chi\in C^\infty(B_r,\R^{N\times 1}\otimes\R^n)$} such that
		\begin{equation*}\label{eq: G-N eq1 Hodge decom 02}
			S(\nabla m,R)=\nabla^\perp\alpha+\chi.
		\end{equation*} 
		Moreover, 
		\begin{equation*}\label{eq: alpha eq 2}
			\Delta\alpha =\pi_n(d\mathbb{P}_R\wedge(dm|0))-\Div^\perp\left[\pi_n(2R\mathbb{P}^{2}(\mathbbm{1}_{n}|0))\right],
		\end{equation*}
		where the linear map $\mathbb{P}_R$ is defined by $\xi\mapsto2R\mathbb{P}(R^{\mathrm{T}}\xi)$.
		
		Similarly, {$\pi_n(d\mathbb{P}_R\wedge(dm|0))\in M_{\loc}^{1,n-2+2\gamma}$} and $ \Div^\perp\left[\pi_n(2R\mathbb{P}^{2}(\mathbbm{1}_{n}|0))\right]\in M_{\loc}^{2,n-2+2\gamma}$. It follows again from Proposition \ref{prop: Holder for Morrey space} {and Proposition \ref{prop:Morrey 2,n-2+2gamma}} that 
		\begin{equation}\label{eq:f alpha L^2zeta estimate 1}
			\begin{aligned}
				\|\nabla\alpha \|_{M^{2\zeta,n-2+\gamma}_*(B_r)}&\lesssim {\|\pi_n(d\mathbb{P}_R\wedge(dm|0))}-\Div^\perp\left[\pi_n(2R\mathbb{P}^{2}(\mathbbm{1}_{n}|0))\right]\|_{M^{1,n-2+\gamma}(B_r)} \\
				&\lesssim \|\nabla R\|_{M^{2,n-2}(B_{r})}\|\nabla m\|_{M^{2,n-2+2\gamma}(B_{r})}+\|\nabla R\|_{M^{1,n-2+\gamma}(B_{r})}\\
				&\lesssim \e\left( \|\nabla R\|_{M^{2,n-2+2\gamma}(B_{r})}+ \|\nabla m\|_{M^{2,n-2+2\gamma}(B_{r})}\right)\\
				&\lesssim \e\left( \|f\|_{L^p(B_{1})}+ 1\right).
			\end{aligned}
		\end{equation}

		Since $R_3,\chi\in C^{\infty}$, $\nabla m,\nabla R \in L^{2\zeta,\infty}(B_r)$. Using standard estimates for harmonic functions and \eqref{eq:f R1 L^2zeta estimate 1}, we infer that for all $0<\rho<r\leq \frac{1}{2}$, there holds 
		$$
		\begin{aligned}
			&\|\nabla R\|_{L^{2\zeta,\infty}(B_\rho)}\lesssim \|\nabla R_3\|_{L^{2\zeta,\infty}(B_\rho)}+ \|\nabla R_1\|_{L^{2\zeta,\infty}(B_\rho)}+ \|\nabla R_2\|_{L^{2\zeta,\infty}(B_\rho)}\\
			\lesssim& \left(\frac{\rho}{r}\right)^{\frac{n}{2\zeta}} \|\nabla R_3\|_{L^{2\zeta,\infty}(B_r)}+ \|\nabla R_1\|_{L^{2\zeta,\infty}(B_r)}+ \|\nabla R_2\|_{L^{\frac{np}{n-p}}(B_r)}\\
			\lesssim&\left(\frac{\rho}{r}\right)^{\frac{n}{2\zeta}} \|\nabla R\|_{L^{2\zeta,\infty}(B_r)}+ \e\left( \|f\|_{L^p(B_1)}+1\right)+ \|f\|_{L^{p}(B_1)}\\
			\lesssim& \left(\frac{\rho}{r}\right)^{\frac{n}{2\zeta}} \|\nabla R\|_{L^{2\zeta,\infty}(B_r)}+ \|f\|_{L^{p}(B_1)}+\e.
		\end{aligned}
		$$Similarly, with \eqref{eq:f alpha L^2zeta estimate 1}, we derive
		$$
		\begin{aligned}
			&\|S(\nabla m,R)\|_{L^{2\zeta,\infty}(B_\rho)}\lesssim \|\nabla \alpha\|_{L^{2\zeta,\infty}(B_\rho)}+ \|\chi\|_{L^{2\zeta,\infty}(B_\rho)}\\
			\lesssim& \left(\frac{\rho}{r}\right)^{\frac{n}{2\zeta}} \|\chi\|_{L^{2\zeta,\infty}(B_r)}+ \|\nabla \alpha\|_{L^{2\zeta,\infty}(B_r)}\\
			\lesssim& \left(\frac{\rho}{r}\right)^{\frac{n}{2\zeta}} \|S(\nabla m,R)\|_{L^{2\zeta,\infty}(B_r)}+ \e\|f\|_{L^{p}(B_1)}+\e.
		\end{aligned}
		$$
		Thus a standard iteration argument gives 
		$$
		\|\nabla R\|_{L^{2\zeta,\infty}(B_r)}+\|S(\nabla m,R)\|_{L^{2\zeta,\infty}(B_r)}\lesssim \|f\|_{L^{p}(B_1)}+\e.
		$$
		Note that $\frac{n}{2\zeta}=n-p$. With \eqref{eq: S estimate 0} we have the estimate
		\begin{equation}\label{eq:f R-m L^2zeta estimate}
			\begin{aligned}
				&\quad\|\nabla m\|_{L^{2\zeta,\infty}(B_r)}+\|\nabla R\|_{L^{2\zeta,\infty}(B_r)}\lesssim  \|f\|_{L^p(B_1)}+\e+r^{n-p}\lesssim  \|f\|_{L^p(B_1)}+1.
			\end{aligned}
		\end{equation}

		If $2\zeta > \frac{np}{n-p}$, then H\"older's inequality implies $\nabla m, \nabla R\in L^{\frac{np}{n-p}}(B_r)$ with the estimate
		$$
		\|\nabla m\|_{L^{\frac{np}{n-p}}(B_r)}+\|\nabla R\|_{L^{\frac{np}{n-p}}(B_r)}\lesssim \|f\|_{L^p(B_1)}+1.
		$$
		
		If $2\zeta\le \frac{np}{n-p}$, then $\nabla m, \nabla R\in L^q_{\loc}$ for any $q\in(2,\frac{4\zeta}{1+\zeta})\subset(2,2\zeta)$. By the definition of $\Omega_{R}$, $\mathbb{P}$, (\ref{eq: S estimate 0}) and H\"older's inequality, we have
		$$
		\Omega_R\cdot\nabla R\in M_{\loc}^{1,n-2+\gamma}\cap L_{\loc}^{\frac{q}{2}},
		$$
		$$
		\mathrm{skew} \left(\nabla m\circ S(\nabla m,R)\right)R \in M_{\loc}^{1,n-2+\gamma}\cap L_{\loc}^{\frac{q}{2}},
		$$
		and
		$$
		{\pi_n(d\mathbb{P}_R\wedge(dm|0))}-\Div^\perp\left[\pi_n(2R\mathbb{P}^{2}(\mathbbm{1}_{n}|0))\right] \in M_{\loc}^{1,n-2+\gamma}\cap L_{\loc}^{\frac{q}{2}}.
		$$
		
		Next, we estimate the $L^{q\zeta}$ norms of $\nabla R_1$ and $\nabla\alpha$. By \eqref{eq:f R1 L^2zeta estimate 1} and \eqref{eq:f alpha L^2zeta estimate 1}, we have 
		\begin{equation}\label{eq:f iteration estimate 1}
			\begin{aligned}
				&\quad\|\Omega_{R}\cdot\nabla R+\mathrm{skew} \left(\nabla m\circ S(\nabla m,R)\right)R)\|_{M^{1,n-2+\gamma}(B_{r})}\\
				&\lesssim \e\left( \|f\|_{L^p(B_1)}+1\right) 
			\end{aligned}
		\end{equation}
		and
		\begin{equation}\label{eq:f iteration estimate 2}
			\begin{aligned}
				&\quad\|\pi_n(d\mathbb{P}_R\wedge(dm|0))-\Div^\perp\left[\pi_n(2R\mathbb{P}^{2}(\mathbbm{1}_{n}|0))\right]\|_{M^{1,n-2+\gamma}(B_r)} \\
				&\lesssim \e\left( \|f\|_{L^p(B_1)}+1\right).
			\end{aligned}
		\end{equation}
		Since $q/2>1$, using H\"older's inequality and \eqref{eq:f R-m L^2zeta estimate}, we have 
		\begin{equation}\label{eq:f iteration estimate 3}
			\begin{aligned}
				&\quad\|\Omega_{R}\cdot\nabla R+\mathrm{skew} \left(\nabla m\circ S(\nabla m,R)\right)R)\|_{L^{\frac{q}{2}}(B_{r})}\\
				&\lesssim\|\Omega_{R}\cdot\nabla R+\mathrm{skew} \left(\nabla m\circ S(\nabla m,R)\right)R)\|_{L^{\frac{2\zeta}{1+\zeta},\infty}(B_{r})}\\
				&\lesssim \|\nabla R\|_{L^2(B_{r})}\|\nabla R\|_{L^{2\zeta,\infty}(B_{r})}+ \| S(\nabla m,R)\|_{L^2(B_{r})}\|\nabla m\|_{L^{2\zeta,\infty}(B_{r})}\\
				&\lesssim \|\nabla R\|_{L^{2\zeta,\infty}(B_{r})}+ \|\nabla m\|_{L^{2\zeta,\infty}(B_{r})}\\
				&\lesssim \|f\|_{L^p(B_1)}+1.
			\end{aligned}
		\end{equation}
		Similarly, using \eqref{eq:f R-m L^2zeta estimate}, we obtain 
		\begin{equation}\label{eq:f iteration estimate 4}
			\begin{aligned}
				&\quad\|{\pi_n(d\mathbb{P}_R\wedge(dm|0))}-\Div^\perp\left[\pi_n(2R\mathbb{P}^{2}(\mathbbm{1}_{n}|0))\right]\|_{L^{\frac{q}{2}}(B_r)} \\
				&\lesssim \|{\pi_n(d\mathbb{P}_R\wedge(dm|0))}-\Div^\perp\left[\pi_n(2R\mathbb{P}^{2}(\mathbbm{1}_{n}|0))\right]\|_{L^{\frac{2\zeta}{1+\zeta},\infty}(B_{r})}\\
				&\lesssim \|\nabla R\|_{L^2(B_{r})}\|\nabla m\|_{L^{2\zeta,\infty}(B_{r})}+ \|\nabla R\|_{L^2(B_{r})}\|\nabla R\|_{L^{2\zeta,\infty}(B_{r})}\\
				&\lesssim \e\left( \|\nabla R\|_{L^{2\zeta,\infty}(B_{r})}+ \|\nabla m\|_{L^{2\zeta,\infty}(B_{r})}\right) \\
				&\lesssim \e\left( \|f\|_{L^p(B_1)}+1\right).
			\end{aligned} 
		\end{equation}
		Applying Proposition \ref{riesz potential M 1,n-beta cap Lp}, \eqref{eq:f iteration estimate 1} and \eqref{eq:f iteration estimate 3}, we infer
		\begin{equation*}\label{eq:f R1 iteration estimate}
			\begin{aligned}
				\|\nabla R_1\|_{L^{q\zeta}(B_r)}&\lesssim \|\Omega_{R}\cdot\nabla R+\mathrm{skew} \left(\nabla m\circ S(\nabla m,R)\right)R)\|_{M^{1,n-2+\gamma}(B_{r})}^{\frac{1}{2-\gamma}}\\
				&\quad\cdot\|\Omega_{R}\cdot\nabla R+\mathrm{skew} \left(\nabla m\circ S(\nabla m,R)\right)R)\|^{1-\frac{1}{2-\gamma}}_{L^{\frac{q}{2}}(B_{r})}\\
				&\lesssim \|f\|_{L^p(B_1)}+1.
			\end{aligned}
		\end{equation*}
		Similarly, using \eqref{eq:f iteration estimate 2} and \eqref{eq:f iteration estimate 4}, we derive
		\begin{equation*}\label{eq:f alpha iteration estimate}
			\begin{aligned}
				\|\nabla\alpha \|_{L^{q\zeta}(B_r)}&\lesssim \|{\pi_n(d\mathbb{P}_R\wedge(dm|0))}-\Div^\perp\left[\pi_n(2R\mathbb{P}^{2}(\mathbbm{1}_{n}|0))\right]\|_{M^{1,n-2+\gamma}(B_r)} ^{\frac{1}{2-2\gamma}}\\
				&\quad\cdot\|{\pi_n(d\mathbb{P}_R\wedge(dm|0))}-\Div^\perp\left[\pi_n(2R\mathbb{P}^{2}(\mathbbm{1}_{n}|0))\right]\|^{1-\frac{1}{2-2\gamma}}_{L^{\frac{q}{2}}(B_{r})}\\
				&\lesssim \e\left( \|f\|_{L^p(B_1)}+1\right).
			\end{aligned}
		\end{equation*}
		Thus we find the following iteration:
		$$
		\nabla m,\nabla R\in L^{q}(B_r)\Longrightarrow \nabla m,\nabla R\in L^{q\zeta}(B_r)
		$$
		together with the estimate 
		\begin{equation*}\label{eq:f R-m q-iteration zeta}
			\|\nabla m\|_{L^{q\zeta}(B_r)}+\|\nabla R\|_{L^{q\zeta}(B_r)}\le C_0\left( \|f\|_{L^p(B_1)}+1\right),
		\end{equation*}
		where $C_0$ is independent of $q$.
		
		Since $\zeta>1$, there exists some $k\ge1$ such that $\zeta^kq\le \frac{np}{n-p}< \zeta^{k+1}q$. After finitely many times iterations, we shall have $\nabla m,\nabla R\in L^{\frac{np}{n-p}}(B_r)$ with the estimate
		\begin{equation*}\label{eq:f R-m p* estimate}
			\|\nabla m\|_{L^{\frac{np}{n-p}}(B_r)}+\|\nabla R\|_{L^{\frac{np}{n-p}}(B_r)}\lesssim \|f\|_{L^p(B_1)}+1.
		\end{equation*}
		
		Now we can estimate the $L^p$ norms of $\Omega_{R}\cdot\nabla R$ and $\mathrm{skew} \left(\nabla m\circ S(\nabla m,R)\right)R$ as follows:
		$$
		\begin{aligned}
			\|\Omega_{R}\cdot\nabla R\|_{L^p(B_r)}&\lesssim \|\nabla R\|_{L^n(B_r)}\|\nabla R\|_{L^{\frac{np}{n-p}}(B_r)}\\
			&\lesssim \|\nabla R\|_{L^{\frac{np}{n-p}}(B_r)}^2
			\lesssim \left( \|f\|_{L^p(B_1)}+1\right)^2 
		\end{aligned}
		$$
		and 
		$$
		\begin{aligned}
			\|\mathrm{skew} \left(\nabla m\circ S(\nabla m,R)\right)R\|_{L^p(B_r)}&\lesssim \|\nabla m\|_{L^n(B_r)}\|S(\nabla m,R)\|_{L^{\frac{np}{n-p}}(B_r)}\\
			&\stackrel{\eqref{eq: S estimate 0}}{\lesssim} \|\nabla m\|_{L^{\frac{np}{n-p}}(B_r)}\left( \|\nabla m\|_{L^{\frac{np}{n-p}}(B_r)}+\|1\|_{L^{\frac{np}{n-p}}(B_r)}\right) \\
			&\lesssim \|\nabla m\|_{L^{\frac{np}{n-p}}(B_r)}^2+r^{\frac{2(n-p)}{p}}
			\lesssim \left( \|f\|_{L^p(B_1)}+1\right)^2.
		\end{aligned}
		$$
		Thus, we have $\Delta R\in L^p$ via \eqref{detailed eq:G-N system 4th eq}. Furthermore, by the usual elliptic regularity theory, \eqref{eq:f R-m M^2,n-2+2gamma estimate} and the above estimates, we have 
		\begin{equation}\label{eq:f R W2,p estimate}
			\|R\|_{W^{2,p}(B_r)}\lesssim \left( \|f\|_{L^p(B_1)}+1\right)^2.
		\end{equation}
		
		By~\cite[Section 6.2]{Gastel-Neff 2022}, the linear operator $L_R:\xi\mapsto\pi_{n}(2R\mathbb{P}^2(R^{\mathrm{T}}(\xi|0)))$ is uniformly positive with
		$$
		\begin{aligned}
			\langle L_R(\xi),\xi\rangle&=\langle\pi_{n}(2R\mathbb{P}^2(R^{\mathrm{T}}(\xi|0))),\xi\rangle=\langle2R\mathbb{P}^2(R^{\mathrm{T}}(\xi|0)),(\xi|0)\rangle\\
			&=\langle2\mathbb{P}(R^{\mathrm{T}}(\xi|0)),\mathbb{P}(R^{\mathrm{T}}(\xi|0))\rangle\ge2\min\{\mu_1,\mu_2,\kappa\}|R^{\mathrm{T}}(\xi|0)|^2=2\min\{\mu_1,\mu_2,\kappa\}|\xi|^2, 
		\end{aligned}
		$$where $\mu_1,\mu_2,\kappa$ is the constant in the definition of operator $\mathbb{P}$. 
		
		Observe that \eqref{detailed eq:G-N system 3rd eq} can be rewritten as an elliptic equation
		\begin{equation*}\label{detailed eq:G-N system 3rd eq rewrite}
			\Div L_{R}(\nabla m)=\Div (\pi_n(2R\mathbb{P}^{2}(R^{\mathrm{T}}(\mathbbm{1}_{n}|0)))).
		\end{equation*}
		Since the coefficients are H\"older continuous, \eqref{eq:f R-m W^2,p estimate} follows from the Calderon-Zygmund theory and \eqref{eq:f R W2,p estimate}.
		\smallskip 
		
		\textbf{Case 2.} $p\ge n$. 
		\smallskip 
		
		In this case, $f\in L^q$ for any $q\in (\frac{n}{2},n)$. Repeating the previous arguments, we conclude that $\nabla m,\nabla R\in L^{\frac{nq}{n-q}}_{\loc}$ with the estimate$$
		\begin{aligned}
			\|\nabla m\|_{L^{\frac{nq}{n-q}}(B_r)}+\|\nabla R\|_{L^{\frac{nq}{n-q}}(B_r)}
			\lesssim \|f\|_{L^p(B_1)}+1.
		\end{aligned}
		$$
		This implies that $R_1,\alpha\in\bigcap_{1<s<\infty}W^{1,s}_{\loc}$ and so $m, R\in W_{\loc}^{2,p}$ with \eqref{eq:f R-m W^2,p estimate} by a similar argument as in Case 1.
	\end{proof}

	
	\bigskip 
	\textbf{Acknowledgements}. C.-Y. Guo and M.-L. Liu are supported by the Young Scientist Program of the Ministry of Science and Technology of China (No.~2021YFA1002200), the NSFC grants (No.~12101362, 12311530037), the Taishan Scholar Project and the NSF of Shandong Province (No.~ZR2022YQ01). C.-L. Xiang is supported by the NFSC grant (No.~12271296) and  the NSF of Hubei province (No.~2024AFA061). The authors would like to thank the anonymous referees for their helpful comments and suggestions that improved our exposition.

	
	

\end{document}